\documentclass[12pt,leqno]{article} 
\usepackage{amsmath,amsfonts}

\newcommand{\bl}{\begin{lemma}}
\newcommand{\bp}{\begin{proposition}}
\newcommand{\bt}{\begin{theorem}}
\newcommand{\bc}{\begin{corollary}}

\newtheorem{proposition}{Proposition}[subsection] 
\newtheorem{lemma}[proposition]{Lemma} 
\newtheorem{corollary}[proposition]{Corollary} 
\newtheorem{theorem}[proposition]{Theorem}

\newcommand{\el}{\end{lemma}}
\newcommand{\ep}{\end{proposition}}
\newcommand{\et}{\end{theorem}}
\newcommand{\ec}{\end{corollary}}

\newtheorem{sublemma}{Lemma}

\newcommand{\prgg}{  \addtocounter{proposition}{1}   
                      \paragraph{\theproposition } }

\newcommand{\per}[1]{\overline{#1}} 

\newcommand{\sfm}[1]{\hbox{{\sf #1}}}

\newcommand{\IH}{I\! \!H}

\newcommand{\TO}{\longrightarrow}
\newcommand{\nt}{\noindent}

\newcommand{\refp}[1]{(\ref{#1})}

\newcommand{\pro}{\noindent {\it Proof. \/}}
\newcommand{\qed}{\hfill $\clubsuit$ \medskip}

\newcommand{\Hom}{\mathop{\rm Hom \, } \nolimits}

\newcommand{\inte}{\mathop{\rm int \, } \nolimits}
\newcommand{\codim}{\mathop{\rm codim \, } \nolimits}
\newcommand{\depth}{\mathop{\rm depth \, } \nolimits}

 \newcommand{\Diff}{\mathop{\rm Diff \, } \nolimits} 

\renewcommand{\c}{\mathfrak{c}}
\newcommand{\I}{\Theta}
\newcommand{\J}{\Lambda}
\newcommand{\C}{\mathfrak{C}}

\newcommand{\rond}{\hbox{\small\hbox{$\circ$}}}

\newcommand{\mf}{M / \mathcal{F}}

\newcommand\phii{{\raise2pt\hbox{$\varphi$}}}  
 \newcommand\phibas{{\raise2pt\hbox{\footnotesize $\varphi$}}}
 \newcommand\Om{\Omega}   
 \newcommand\om{\omega} 
 
\catcode`\‡=\active\def‡{\'a}        
\catcode`\Ž=\active\defŽ{\'e}        
\catcode`\'=\active\def'{\'{\i}}     
\catcode`\—=\active\def—{\'o}        
\catcode`\œ=\active\defœ{\'u}        
\catcode`\ƒ=\active\defƒ{\'E}        
\catcode`\ˆ=\active\defˆ{\`a}        
\catcode`\=\active\def{\`e}        
\catcode`\"=\active\def"{\`{\i}}     
\catcode`\˜=\active\def˜{\`o}        
\catcode`\=\active\def{\`u}        
\catcode`\Ë=\active\defË{\`A}        
\catcode`\‹=\active\def‹{\~a}        
\catcode`\–=\active\def–{\~n}        
\catcode`\›=\active\def›{\~o}        
\catcode`\Ì=\active\defÌ{\~A}        
\catcode`\"=\active\def"{\~N}        
\catcode`\Í=\active\defÍ{\~O}        
\catcode`\‰=\active\def‰{\^a}        
\catcode`\=\active\def{\^e}        
\catcode`\"=\active\def"{\^{\i}}     
\catcode`\™=\active\def™{\^o}        
\catcode`\ž=\active\defž{\^u}        
\catcode`\°=\active\def°{\infty}      

 \renewcommand\C{\mathbb C}

 \newcommand\R{\mathbb R}

 \renewcommand\S{\mathbb S}
 \newcommand\T{\mathbb T}
 
 \newcommand\D{\mathbb D}
 \newcommand\Z{\mathbb Z}

 \newcommand{\sbat}{{\S}^{^1}}
 
 \newcommand{\bi}[2]{{#1}^{^{#2}}}
 \newcommand{\hiru}[3]{{#1}^{^{#2}}{\left( #3 \right)}}
 \newcommand{\lau}[4]{{#1}^{^{#2}}_{_{#3}}{\left( #4 \right)}}
 
\title{\Large\bf  {The  BIC of a conical fibration.}}

\author{
Martintxo Saralegi-Aranguren\thanks{UPRES-EA 2462 GŽomŽtrie-Algbre. 
FacultŽ Jean Perrin.
UniversitŽ d'Artois.   Rue Jean Souvraz SP 18.   62 307 Lens Cedex - 
France.   
{\sl saralegi@euler.univ-artois.fr}. Partially supported by the 
"Programa 
de 
perfeccionamiento y movilidad de personal investigador" of the 
"Eusko Jaularitza-Gobierno del Pa's Vasco" and "Projet d'Action IntegrŽe 
Picasso" n. 02679SL}
\\ {\small UniversitŽ d'Artois }
\and   
Robert Wolak\thanks{Instytut Matematyki. Uniwersytet Jagiellonski. 
Wl. 
Reymonta4, 30 059
Krakow - Poland.  
{\sl wolak@im.uj.edu.pl}.  Partially supported by the KBN grant 5
P03A 040 20}
\\ {\small Uniwersytet Jagiellonski}
}
\date{\today}
\begin{document}  
\maketitle

\begin{abstract}  

In the paper we introduce the notions of a singular fibration and a
singular Seifert fibration. These notions are natural generalizations of
the notion of a locally trivial fibration to the category of stratified
pseudomanifolds. For singular foliations defined by such fibrations we
prove a de Rham type theorem for the basic intersection cohomology
introduced the authors in a recent paper. One  of important examples of
such a structure is the natural projection onto the leaf space for the
singular Riemannian foliation defined by an action of a compact Lie
group on a compact smooth manifold.
\end{abstract}

The  failure  of the PoincarŽ duality for the homology and cohomology 
of some singular spaces
led Goresky and MacPherson to introduce a new homology theory called 
intersection homology which
took into account the properties of the singularities of the 
considered space (cf. \cite{GM1}). 
This homology is defined for stratified pseudomanifolds.
The initial idea was generalized in several ways. 
The theories of simplicial and singular homologies were developed as 
well as weaker notions of the perversity were proposed (cf. 
\cite{K,M}). Several versions of the PoincarŽ duality were proved 
taking into account the notion of dual perversities. Finally, the deRham 
intersection cohomology was defined by Goresky and MacPherson 
for Thom-Mather stratified spaces. The first written reference is the 
paper by J.L. Brylinski, cf. \cite{Bry}. This led to the search for 
the "de Rham - type"  theorem  for the intersection cohomology.

The first author in his thesis and several subsequent publications, 
(cf. 
\cite{S}) presents the de Rham intersection cohomology of stratified 
spaces using essentially 
the Verona resolution of singularities (cf. \cite{V}). 
The main idea is the introduction of the complex (sheaf) of liftable 
intersection 
forms (cf. \cite{S}).

In \cite{SW1} the authors introduced the basic (de Rham) intersection 
cohomology (BIC)
for singular Riemannian foliations  using the fact that such 
foliations define
a stratification of the manifold (cf. \cite{Mol,W1}) 
respect to the leaves of the foliation. In fact, the  BIC can be defined 
for a larger class
of foliations which we call conical foliations.

In the present paper  we study BIC for a particular class of
conical foliations which are defined by the so called singular 
fibration.
In this case the space of leaves is a stratified pseudomanifold and 
the BIC of this foliation is the intersection cohomology of the space 
of leaves. This result is a version   of the well-known de Rham 
theorem  for the basic de Rham intersection cohomology of a foliation 
defined by a generalized fibration.

Such singular fibrations can be met in several natural situations, 
for example the foliation defined by an action of a compact Lie group 
is defined by a singular fibration (cf. \cite{D}).

\section{Conical foliations.}

\subsection {\bf Singular foliations.}

A {\em regular  foliation}   on a manifold $M$ is a partition 
${\mathcal  F}$ of $M$  by connected immersed submanifolds, called 
{\em leaves}, with the following local model: 
$$(\R^m,{\mathcal H})$$ 

\noindent
where leaves of $\cal H$ are defined by $\{ dx_1 = \cdots = dx_p = 
0\}$. We shall say that $(\R^m,{\mathcal H})$ is a {\em simple 
foliation}. Notice that the leaves have the same dimension.

A {\em singular foliation}  on a manifold $M$ is a partition 
${\mathcal  F}$ of $M$  by connected immersed submanifolds, called 
{\em leaves}, with the following local model:
$$
(\R^{m-n-1} \times \R^{n+1},{\mathcal H} \times {\mathcal K}) $$ 
\noindent
where $(\R^{m-n-1},{\mathcal H})$ is  a simple foliation and 
$(\R^{n+1},{\mathcal K})$ is  a singular foliation having the origin 
as a leaf. When $n=-1$ we just have a regular foliation. Notice that 
the leaves may have different dimensions. This local model is exactly 
the local model of a foliation of Sussmann \cite{Su} and Stefan 
\cite{St}; so, they are singular foliations in our sense 

  Classifying the points of $M$ following the dimension of the leaves 
one gets a stratification $ {\sfm{S}}_{\mathcal   F}$ of $M$ whose 
elements are called {\em strata}. The foliation is regular  when this 
stratification has just one stratum
$\{ M\}$. 

A smooth map $f \colon (M,{\mathcal  F}) \rightarrow  (M',{\mathcal  
F'})$  between singular foliated manifolds is {\em foliated} if it 
sends the leaves of ${\mathcal  F}$ into the leaves of ${\mathcal  
F'}$. When $f$ is an embedding it preserves the dimension of the 
leaves and therefore it sends the strata of $M$ into the strata of 
$M'$. The group of foliated diffeomorphisms preserving the foliation 
will be denoted by $\Diff (M, {\mathcal  F})$.  Examples, more 
properties and the singular version of the Frobenius theorem the 
reader can find in \cite{SW1,BA,St,Su,VA}.

\prgg {\bf Examples}.

(a) In any open subset $U \subset M$ we have the singular foliation
${\mathcal  F}_U = $\{connected components of  $L \cap U \ / \ L \in 
{\mathcal  F}$ \}. The associated stratification is $
{\sfm{S}}_{{\mathcal  F}_U} = \{ \hbox{connected components of } S
\cap U \ / \ S\in  {\sfm{S}}_{\mathcal   F}\}. $

\medskip

(b) We shall say that a foliated embedding of the form 
$$ \phii \colon (\R^{m-n-1} \times \R^{n+1},{\mathcal H} \times 
{\mathcal K}) \longrightarrow (U,{\mathcal  F}_U), $$  is a {\em 
foliated chart}  modelled on $(\R^{m-n-1} \times \R^{n+1},{\mathcal 
H} \times {\mathcal K})$.

\medskip

(c) Consider $(N,{\mathcal K})$ a connected regular foliated 
manifold. In the product $N \times M$ we have the singular foliation 
${\mathcal K} \times
{\mathcal  F} = \{ L_1 \times L_2 \ / \ L_1 \in {\mathcal K},  L_2 
\in {\mathcal  F}\}$.  The associated stratification is $ {\sfm{S}}_{{\mathcal K} \times {\mathcal  F}} = \{  N \times  S \ / \  S \in 
{\sfm{S}}_{\mathcal   F}\}.$

\medskip

 (d) Consider $\S^n$ the sphere of dimension $n$
 endowed with a singular foliation ${\mathcal G}$ without 
 $0$-dimensional leaves. Identify 
the disk $D^{n+1}$ with the cone $c\S^n = \S^n \times [0,1[  \ / \ \S^n 
\times \{ 0\}$ by the map  $ x \mapsto [x/||x||,||x||]$ where  
$[u,t]$ is a generic element of the cone. We shall consider on 
$D^{n+1}$  the singular foliation 
$$ c {\mathcal G} =  \{ F \times \{ t \}  \ / \  F \in {\mathcal G}, 
t \in ]0,1[ \} \cup \{ \vartheta \},$$  
where $\vartheta$ is the vertex $[u,0]$ of the cone. The induced 
stratification is 
$$ {\sfm{S}}_{c{\mathcal G}} = \{ S \times ]0,1[  \ / \  S \in {\cal 
S}_{{\mathcal G}}\} \cup \{\vartheta\}, $$ 
since $\mathcal{G}$ does not possesses $0$-dimensional leaves.
For technical reasons we allow $n=-1$, in this case 
$\S^n = \emptyset$ and $c\S^n  = \{ \vartheta \}$.

\smallskip

Unless otherwise stated, if $(\S^n,{\mathcal G})$ is a singular 
foliation, 
we shall consider on $D^{n+1} = c \S^n$  the foliation $c{\mathcal  
G}$.

\subsection {\bf Conical foliations.}

The strata of this kind of foliations are 
not necessarily manifolds  and their relative position can be very 
wild. 
Consider $(\R,{\mathcal  F})$ where ${\mathcal  F}$ is given by a 
vector field $f
\frac{\partial}{\partial t}$; there are two kind of strata. The 
connected
components of $f^{-1}(\R-\{ 0\})$ and these of $f^{-1}(0)$. 
In other words, any connected closed subset of $ S$ can be a stratum. 
In order  to support a intersection cohomology structure, the 
stratification
${\mathcal S}_{{\mathcal  F}}$ asks for a certain amount of 
regularity and 
conicalicity (see \cite{GM1} for the case of  stratified 
pseudomanifolds). This leads us to introduce the following definition.

A singular foliation $(M,{\mathcal  F})$ is said to be a {\em conical 
foliation} if  any point $x \in M$  possesses a foliated chart 
$(U,\phii)$ modelled on
$$(\R^{m-n-1} \times c\S^n,{\mathcal H} \times c{\mathcal G}),$$
where
$(\S^n,{\mathcal G})$ is a conical foliation without  0-dimensional 
leaves. We shall say that $(U,\phii)$ is a {\em conical chart} of 
$x$ and that $(\S^n,{\mathcal G})$ is the {\em link} of $x$.  
Notice that, if $S$ is the stratum containing $x$ then 
$\phii^{-1}(S \cap  U) = \R^{m-n-1} \times \{ \vartheta \}$. This 
definition is made by induction on the dimension of $M$.

 Notice that each stratum is an embedded submanifold of $M$. The leaf 
(resp. stratum) of ${\mathcal  F}$ (resp. ${\sfm{S}}_{\mathcal   F}$) 
containing $x$ is  sent by $\phii$ to the leaf of ${\mathcal H}$ 
containing $0$ (resp. $\R^{m-n-1}$). Although a point $x$ may have 
several charts the integer $n+1$ is an invariant: it is the 
codimension of the stratum containing $x$. This integer cannot to be 
1 since the conical foliation $(\S^n,{\mathcal G})$  has not 
0-dimensional leaves.

We also use the notion of conical foliation in a manifold with 
boundary. In this case, the points of the boundary have conical 
charts modelled on $(\R^{m-n-2} \times [0,1[ \times c\S^n,{\mathcal H} 
\times {\mathcal I} \times c{\mathcal G})$ where ${\mathcal I}$ is 
the foliation by points of $[0,1[$. The boundary $(\partial M, 
{\mathcal  F}_{\partial M})$ is also a conical foliated manifold.

The above local description implies some important facts about the 
stratification ${\cal 
S}_{\mathcal  F}$. Notice for example that the
family of strata is finite in the compact case and locally finite 
in the general case.
The closure of a stratum $S\in {\cal 
S}_{\mathcal   F}$ is a union of strata. Put $ S_1 \preceq S_2$ if $S_1, S_2 \in {\cal 
S}_{\mathcal   F}$ and
$S_1 \subset \overline{S_2}$. This relation is an order relation and 
therefore
$({\sfm{S}}_{\mathcal   F},\preceq)$ is a poset. 

The  $\depth$  of ${\sfm{S}}_{\mathcal   F}$, written $ \depth  {\cal 
S}_{\mathcal   F}$, is defined to be the largest $i$ for which there 
exists a chain of strata $S_0  \prec S_1 \prec \cdots \prec S_i$. So, 
$ \depth  {\sfm{S}}_{\mathcal  F} = 0$ if and only if the foliation 
${\mathcal  F}$ is regular. We also have $ \depth  {\cal 
S}_{{\mathcal  F}_U} \leq  \depth   {\sfm{S}}_{{\mathcal  F}}$ for any 
open subset $U \subset M$ and  $ \depth  {\sfm{S}}_{{\mathcal G}} =  \depth  {\cal 
S}_{{\mathcal H} \times {\mathcal G}} < \depth  {\sfm{S}}_{{\mathcal 
H} \times c{{\mathcal G}}}$ (cf. 1.1.1).

The minimal strata are exactly the closed strata. An inductive 
argument shows that the maximal strata are the strata of dimension 
$m$. They are called {\em  regular strata} and the others {\em 
singular strata}.   The union of singular strata is written $\Sigma_{\mathcal  F} $. Since the codimension of singular strata is at 
least 2, then only one regular $R$ strata appears, which is an open 
dense subset.
\bigskip

 For the definition of perverse forms we exploit the compatibility 
between the foliated charts  given in the following Lemma whose proof 
can be found in \cite{SW1}.

\begin{lemma}
\label{lema} 
Let $(U_1,\phii_1)$, $(U_2,\phii_2)$ be two foliated charts of a point 
$x$ 
of $M$  with $U_1 \subset U_2$. There exists  an unique foliated 
embedding $$
\phii_{1,2} \colon (\R^{m-n-1} \times  \S^n \times [0,1[,{\mathcal H}_1 
\times 
{\mathcal G}_1 \times {\mathcal I} ) \rightarrow (\R^{m-n-1} \times  
\S^n \times [0,1[,{\mathcal H}_2 \times {\mathcal G}_2   \times 
{\mathcal I})$$
making the following diagram commutative 
$$ \begin{picture}(170,60)(00,-5)

\put(0,0){\makebox(0,0){$ \R^{m-n-1} \times  c \S^n $}}
\put(180,0){\makebox(0,0){$  \R^{m-n-1} \times  c \S^n $}}
\put(180,50){\makebox(0,0){$  \R^{m-n-1} \times   \S^n \times [0,1[$}}
\put(0,50){\makebox(0,0){$ \R^{m-n-1} \times   \S^n \times [0,1[$}}
\put(0,40){\vector(0,-1){26}} 
\put(180,40){\vector(0,-1){26}}
\put(55,0){\vector(1,0){66}}
\put(75,50){\vector(1,0){46}}
\put(88,12){\makebox(0,0){$\phii_2^{-1} \rond \phii_1$}}
\put(98,58){\makebox(0,0){$\phii_{1,2}$}}
\put(-10,27){\makebox(0,0){$P$}}
\put(190,27){\makebox(0,0){$P$}}
\end{picture} 
$$ 
where the smooth map $P$ is defined by $P(u,\theta,t) = (u,[ 
\theta,t])$.
\end{lemma}

This result also implies that the link $(\S^n,\mathcal{G})$ of a 
point $x$ is the 
same for any point of the stratum $S$ containing $x$. For this reason we shall say that 
$(\S^n,\mathcal{G})$ is {\em the link of the stratum} $S$.
\subsection{\bf Partial blow up.}

A useful tool for inductive 
arguments is the desingularisation. Fix $(M,{\mathcal  F})$  a 
conical foliated manifold with strictly positive depth. A {\em 
partial blow up} is  a conical foliated manifold with boundary 
$(\widehat{M},\widehat{{\mathcal  F}})$ and a foliated smooth map 
${\mathcal L}_M \colon  (\widehat{M},\widehat{{\mathcal  F}}) 
\longrightarrow (M,{\mathcal  F})$ such that:

\begin{itemize}
\item[1)] The restriction ${\mathcal L} \colon \widehat{M} - 
\partial\widehat{M} 
\longrightarrow  M - S_{_{min}}$ is a diffeomorphism, where 
$S_{_{min}}$ is the union of closed (minimal) strata.

\item[2)] For any point $x \in S_{_{min}}$ there exists a commutative 
diagram

$$ \begin{picture}(120,60)(00,-8)
\put(0,0){\makebox(0,0){$ \R^{m-n-1} \times c \S^n$}}
\put(120,0){\makebox(0,0){$M$}}
\put(120,50){\makebox(0,0){$\widehat{M}$}}
\put(0,50){\makebox(0,0){$ \R^{m-n-1} \times  \S^n \times [0,1[$}}
\put(0,40){\vector(0,-1){26}} 
\put(120,40){\vector(0,-1){26}} 
\put(50,0){\vector(1,0){56}}
\put(70,50){\vector(1,0){36}}
\put(78,8){\makebox(0,0){$\phii$}} 
\put(88,58){\makebox(0,0){$\widehat\phii$}}
\put(-10,27){\makebox(0,0){$P$}} 
\put(135,27){\makebox(0,0){${\mathcal L}$}} 
\end{picture} 
$$ 
\noindent
where $(U,\phii)$ is a conical chart of $x$ and $\phii \colon ( 
\R^{m-n-1} \times  \S^n \times [0,1[ , {\mathcal H} \times {\mathcal 
G} \times {\mathcal I})\rightarrow (\widehat{M},\widehat{{\mathcal  
F}})$ is a foliated embedding.
\end{itemize}

\medskip

\noindent Notice that each restriction ${\mathcal  L} \colon 
{\mathcal  L}^{-1}(S) \rightarrow S$, where $S$ is a closed stratum, 
is a fiber bundle possessing a foliated 
atlas whose structural group is $\Diff ( \S^n,{\mathcal G})$ relatively 
to a conical foliation $( \S^n,{\mathcal G})$.

\bigskip

\begin{proposition}
The  partial blow up always exists and is unique. It also verifies 
\begin{itemize}
    \item[3)]  $ \depth {\sfm{S}}_{\widehat{{\mathcal  F}}} = 
 \depth {\sfm{S}}_{{\mathcal  F}}-1$.

\noindent \item[4)] $(\partial {\widehat{M}}, \widehat{{\mathcal  
F}}_{\partial{\widehat{M}}})$ is a conical foliated manifold.
\end{itemize}
\end{proposition}

The proof of Proposition 1.3.1 is presented in \cite{SW1}.

\section{Basic Intersection cohomology.}

There are two ways to define perverse forms, the first one uses a 
system of 
tubular neighborhoods (cf. \cite{Bry}) and the second one uses a 
global blow up (\cite{S}). In the study of conical foliations we have 
decided to introduce intersection forms in an intermediate way: 
iterating the local blow up $M$ we obtain  a manifold  
$\widetilde{M}$ with borders where the foliation become regular. A 
similar procedure has been developed in \cite{S}. Some of the  
differential forms on the regular part of $M$  can be extended to 
$\widetilde{M}$, these are the perverse forms. We present this  
notion in an intrinsic way without constructing $\widetilde{M}$. 

\smallskip

We are going to deal with differential forms on  a product $N \times 
[0,1[^p$, where $N$ is a manifold: they are restrictions of 
differential forms defined on $N \times ]-1,1[^p$.

\subsection{\bf Perverse forms.} 

The differential complex 
${\bf\Pi}^{\ast}({M \times [0,1[^p})$ of perverse forms is introduced 
by induction on the $ \depth   {\sfm{S}}_{{\mathcal  F}}$. When this 
depth is  0 then ${\bf\Pi}^{\ast}({M \times [0,1[^p}) = 
\hiru{\Omega}{\ast}{M \times [0,1[^p}$. In the general case we shall 
put $\omega \in {\bf\Pi}^{\ast}({M \times [0,1[^p})$ if $\omega$ is a 
differential form on $ \hiru{\Om}{\ast}{(M - \Sigma_{\mathcal  F} 
\times [0,1[^p})$ such that any $x \in M$ possesses a conical chart  
$(U,\phii)$ such that 

$$
(\phii \times \hbox{identity}_{[0,1[^p})^*\omega  \ \ 
\hbox{ extends to }\ \ \
\omega_\phii \in {\bf\Pi}^{\ast}({ \S^{m-n-1} \times  \S^n   \times 
[0,1[^{p+1}}).
$$

In fact, ${\bf\Pi}^{\ast}({M \times [0,1[^p})$ is a differential 
graduated 
commutative algebra (dgca in short) since $(\omega +\eta)_\phii = 
\omega_\phii + \eta_\phii$, $(\omega \wedge\eta)_\phii = \omega_\phii 
\wedge \eta_\phii$ and $(d\omega)_\phii = d\omega_\phii$.

Notice that the notion of perverse form depends on the foliation 
${\mathcal  F}$ trough the stratification ${\sfm{S}}_{{\mathcal  F}}$.

\medskip

\prgg {\bf Properties} (see \cite{SW1}).

(a) Let $(N,{\mathcal H})$ be a regular foliated manifold. The 
partial blow 
up of $N \times c \S^n$ is $N \times  \S^n \times [0,1[$, with 
$(x,\theta,t)
\mapsto (x , [\theta,t])$. Here the factor $[0,1[$ appears. Further 
desingularisation produces extra $[0,1[$ factors.

\medskip

  (b) Consider ${\mathcal  L}_M \colon \widehat{M} \rightarrow M$  
the partial blow up of a conical foliated manifold. Any perverse form 
$\omega$ on $M$ defines a perverse form $\widehat\omega$ on 
$\widehat{M}$ extending ${\mathcal  L}_M^*\omega$. This map $\omega 
\mapsto \widehat\omega$ defines an isomorphism between 
${\bf\Pi}^{\ast}({M} )$ and  
$\Pi^{^*}(\widehat{M})$.

\medskip


\medskip

(c) Any perverse form and any conical chart verify the equation (1).

\medskip

(d) There are smooth functions on $M -\Sigma_{\mathcal  F}$ which 
are not 
perverse. Any differential form $\omega$ of $M$ is perverse.

\medskip

(e) Consider $\{ U , V \}$ an open covering of $M$. There exist a 
subordinated partition of the unity made up with smooth functions 
defined on $M$. The Mayer-Vietoris short sequence
$$0 \rightarrow {\bf\Pi}^{\ast}({M}) \rightarrow 
{\bf\Pi}^{\ast}({U}) \oplus {\bf\Pi}^{\ast}({V}) \rightarrow 
{\bf\Pi}^{\ast}({U\cap V}) \rightarrow 0, $$
where the map are defined by $\omega \mapsto (\omega,\omega)$ and 
$(\alpha,\beta) \mapsto \alpha - \beta$, is exact.





\subsection{\bf Perverse degree.} 

The amount of transversality of a perverse form $\omega \in 
{\bf\Pi}^{\ast}({M})$ with respect a singular stratum $S$ is measured 
by the perverse degree. Given a point $x$ of a singular stratum $S$, 
we define the {\em perverse degree} $||\omega||_x$ as the smallest 
integer $k$ verifying:
$$
i_{\zeta_{0}} \cdots  i_{\zeta_{k}}\omega_\phii \equiv 0 
$$
for each conical chart $(U,\phii)$ of $x$, with $\phii^{-1}(x) = (a, 
\vartheta )$, and each family of vector fields $  \{ 
\zeta_{i}\}_{i=0}^{k} $ on
$\{ a \}\times ( \S^n - \Sigma_{\mathcal G}) \times \{ 0 \}$. Here 
$i_\cdot$ 
denotes the interior product. When the form $\omega_\phii$ vanishes on 
$]a-\epsilon,a+\epsilon[ \times ( \S^n - \Sigma_{\mathcal G}) \times 
\{ 0 \}$ then we shall write $|| \omega ||_x =- \infty$. We define 
the {\em perverse degree} 
$||\omega||_S$ by
$$
||\omega||_S =  \sup \{ || \omega ||_x \ / \ x \in S \}.
$$

\prgg {\bf Properties}. 

We recall some basic properties of the 
perverse form, cf. \cite{SW1}.

(a) A perverse form $\omega$ can be extended to $\widetilde{M}$. Its 
perverse degree is in fact the degree of $\omega$ on the added part 
when passing from $M$ to $\widetilde{M}$.

\medskip


\medskip

(b) The perversity conditions (2) and (3) do not depend on the choice 
of the conical chart.

\medskip

(c) For two perverse forms $\omega$ and $\eta$ and a singular stratum 
$S$ we have: $$ ||\omega + \eta||_S \leq \max \{ ||\omega||_S , ||  
\eta||_S\} \ , 
||\omega \wedge \eta||_S \leq  ||\omega||_S +||  \eta||_S .$$

\medskip

(d) The perverse degree of a perverse function is 0 or $- \infty$. 
The perverse degree of a  differential form $\omega$ of $M$ is 0 or 
$- \infty$.

\medskip

 (e) Let $(N,{\mathcal H})$ be a regular foliation and $a_0 \in N$ a 
basis point. Put $pr \colon  N \times M \rightarrow M$ the canonical 
projection and $\iota \colon M \rightarrow N \times M$ the inclusion 
defined by $\iota (x) = (a_0,x)$. For any $\omega \in 
{\bf\Pi}^{\ast}({M})$ and $\eta\in{\bf\Pi}^{\ast}({N\times M})$ we 
have 
$$ ||pr^* \omega||_{N \times S} \leq ||\omega||_S \ \hbox{ and } \ 
||\iota^* \omega||_S \leq ||\omega||_{N \times S} $$
for each singular stratum of $M$.

\medskip

(f) Consider $\omega \in {\bf\Pi}^{\ast}({c \S^n})$ a perverse form. 
Its perverse degree relatively to the vertex is:
$$ ||\omega||_{\{ \vartheta \} } = \left\{ \begin{array}{ll}
- \infty & \hbox{ if } \omega_\phii \equiv 0 \hbox{ on }  \S^n \times \{
0 \} \\ \deg \omega & \hbox{ if not,}
\end{array}
\right.$$
where $(U,\phii)$ is any conical chart of $\{ \vartheta \}$.

\medskip



(g) A perverse form with $||\omega||_S \leq 0$ and $ || d\omega ||_S 
\leq 0$ induces a differential form $\omega_S$ on $S$. When this 
happens for each 
stratum we conclude that $\omega \equiv \{ \omega_S\}$ is a Verona's 
{\em
controlled form} (cf. \cite{V}).

\subsection{\bf Basic cohomology.}

Consider $(M, {\mathcal  F})$  a 
foliated manifold. A differential form $\omega \in 
\hiru{\Om}{\ast}{M}$ is {\em basic} if
$$
i_X \omega = i_X d \omega = 0,
$$
for each vector field $X$ on $M$ tangent to the foliation. By 
$\hiru{\Om}{\ast}{\mf}$ we denote the complex of basic forms. 
Since the sum and the product of basic forms are still basic forms 
then the complex of basic forms is a dgca. Its cohomology 
$\hiru{H}{*}{\mf}$ is the {\em basic cohomology} of $(M,{\cal 
F})$, which is a dgca. If $F\subset M$ is a saturated closed subset 
we shall write 
$\hiru{\Om}{\ast}{(M,F)/{\mathcal  F}}$ the subcomplex of basic 
forms of $M$ vanishing on $F$. Its cohomology will be denoted by 
$H^{\ast}{((M,F)/{\mathcal  F})}$

\prgg {\bf Properties}.

(a) A smooth foliated map $f \colon  (M_1,\mathcal{F}_1) \to 
(M_2,\mathcal{F}_2)$  
induces  the dgca operator $f^* \colon \hiru{\Om}{*}{M_2  / 
\mathcal{F}_2}$ $\to
\hiru{\Om}{*}{M_1 / \mathcal{F}_1}$.

\medskip

(b) An open covering $\{ U , V \}$  of $M$ by saturated open 
subsets is a {\em basic covering}. 
when there exists a 
subordinated partition of the unity made up of basic  functions. They 
may or may not exist.
Then the Mayer-Vietoris short sequence

$$
0 \to \hiru{\Om}{*}{\mf} \to \hiru{\Om}{*}{U/ \mathcal{F}_U} \oplus 
\hiru{\Om}{*}{V/ \mathcal{F}_V} \to
\hiru{\Om}{*}{(U\cap V)/ \mathcal{F}_{U\cap V}} \to 0,
$$
where the map are defined by $\om \mapsto (\omega,\omega)$ and 
$(\alpha,\beta) \mapsto \alpha - \beta$, is exact. This result is not 
longer true for 
more general coverings


\subsection{\bf Basic intersection cohomology.}

Consider $(M, {\mathcal  F})$  a conical foliated manifold.  A {\em 
perversity} is a map $\per{p} \colon {\sfm{S}}_{\mathcal   F} 
\rightarrow \Z  $. There are two particular perversities: the {\em 
zero perversity} $\per{0}$ and the {\em top perversity} $\per{t}$ 
defined by $\per{0}(S) =0$ and $\per{t}(S) = \codim S-2$. 
Associated to a smooth foliated map $f \colon (M',{\mathcal  F}') 
\rightarrow (M, {\mathcal  F}) $ (resp. $f \colon (M,{\mathcal  F}) 
\rightarrow (M', {\mathcal  F}') )$ there exists a perversity on 
$(M',{\mathcal  F}') $, still written $\per{p}$, defined by 
$\per{p}(S') = \per{p}(S)$ where $S' \in {\sfm{S}}_{{\mathcal  F}}$ 
and $S \in {\sfm{S}}_{{\mathcal  F}}$ with $f(S') \subset S$ (resp. $S 
= f^{-1}(S')$).

The basic intersection cohomology appears when one considers basic 
forms whose 
perverse degree is controlled by a perversity. We shall write
$$\lau{\Omega}{\ast}{\per{p}}{\mf}  = \{ \omega \in 
\hiru{\Om}{\ast}{(M - 
\Sigma_{\mathcal  F}) / {\mathcal  F}_{M - \Sigma_{\mathcal  F}}} \ / 
\ max (||\omega||_S ,||d\omega||_S )\leq \per{p}(S) \ \ \forall S \in 
{\sfm{S}}_{{\mathcal  F}} \}$$
the complex of basic forms whose perverse degree (and that of the 
their derivative) is bounded by the perversity $\per{p}$. It is a 
differential 
complex, but it is not an algebra, in fact the wedge product acts in 
this way:

$$ \wedge \colon \lau{\Omega}{i}{\per{p}}{\mf} \times 
\lau{\Omega}{j}{\per{q}}{\mf} \rightarrow 
\lau{\Omega}{i+j}{\per{p} + \per{q}}{\mf}$$

\noindent
The cohomology $\lau{\IH}{\ast}{\per{p}}{M/{\cal F}}$ of this complex is 
the {\em basic intersection cohomology} of $M$, or BIC for short, 
relatively to the perversity $\per{p}$.

\prgg {\bf Properties}.

(a) A foliated embedding$f \colon  (M_1,\mathcal{F}_1) \to 
(M_2,\mathcal{F}_2)$ 
induces  the dgca operator $f^* \colon \lau{\Om}{*}{\per{p}}{M_2 / 
\mathcal{F}_2} \to 
\lau{\Om}{*}{\per{p}}{M_1/ \mathcal{F}_2}$.

\medskip



(b) Consider $\{ U , V \}$ a basic covering of $M$.
Then the Mayer-Vietoris short sequence
$$
0 \to \lau{\Om}{*}{\per{p}}{\mf} \to 
\lau{\Om}{*}{\per{p}}{U/\mathcal{F}_U} \oplus 
\lau{\Om}{*}{\per{p}}{V/\mathcal{F}_V}  \to
\lau{\Om}{*}{\per{p}}{(	U \cap V)/\mathcal{F}_{U \cap V}}  \to 0,
$$
where the map are defined by $\om \mapsto (\omega,\omega)$ and 
$(\alpha,\beta) \mapsto \alpha - \beta$, is exact. This result is not 
longer true for 
more general coverings.   

\bigskip

The basic intersection cohomology coincides with the basic cohomology 
when the foliation ${\mathcal  F}$ is regular. But it also 
generalizes the intersection cohomology of Goresky-MacPherson (cf. 
\cite{GM1}) when the leaf space $B$ lies in the right category, that 
of stratified pseudomanifolds (cf. Theorem \ref{Compacto}).

The intersection basic cohomology, as the basic cohomology, are not 
easily computed. Nevertheless, the typical calculations for 
intersection basic cohomology are the following.

\begin{proposition}
\label{local4}
Let $( \R^{k} ,{\mathcal H})$ be a simple foliation and let 
$(M,{\mathcal  F})$ be  a conical foliation. For any perversity 
$\per{p}$ the projection $pr \colon \R^k \times (M - 
\Sigma_{{\mathcal  F}}) \rightarrow (M - \Sigma_{{\mathcal  F}})$ 
induces an isomorphism $$\lau{\IH}{\ast}{\per{p}}{\R^k \times M/ 
{\mathcal H} \times {\mathcal  F}} \cong 
\lau{\IH}{\ast}{\per{p}}{\mf}.$$
\end{proposition}

\begin{proposition}
\label{local5}
Let ${\mathcal G}$ be a conical foliation on the sphere $\S^n$. For 
any perversity $\per{p}$ the projection $pr \colon ( \S^n - 
\Sigma_{{\mathcal G}}) \times ]0,1[ \rightarrow ( \S^n - 
\Sigma_{{\mathcal G}})$ induces the isomorphism
$$\lau{\IH}{i}{\per{p}}{c S^k / c{\mathcal G}} =\left\{
\begin{array}{cl}
\lau{\IH}{i}{\per{p}}{ S^k / {\mathcal G)} }& \hbox{if } i \leq 
\per{p}(\{ 
\vartheta \}) \\
0 & \hbox{if } i > \per{p}(\{ \vartheta \}). 
\end{array}
\right.
$$
\end{proposition}

This result shows that the basic intersection cohomology is not completely 
determinated 
by the cohomology of $M$. 

\section{Stratified pseudomanifolds.}

In some cases, mainly when the leaves are compact, the orbit space 
$B$ of a conical foliated manifold $(M,\mathcal{F})$ has a 
nice topological structure, that of stratified pseudomanifold. 
This notion  has been introduced in 
\cite{GM1}. In the smooth context, a {\em stratified pseudomanifold} 
$X$
is  given by
\begin{itemize}
    \item a paracompact space $X$ 
    \item  a locally finite partition ${\sfm{S}}_{X}$, called  
    {\em stratification}, made up of connected smooth 
    manifolds, called {\em strata}. 
    \end{itemize}
The local structure must be conical, that is, the local model is $\R^p \times 
cL$, where $L$ is a smaller compact stratified pseudomanifold. More 
exactly, each point of $X$ has a neighborhood which is homeomorphic to 
the local model, and this homeomorphism sends diffeomorphically a 
stratum into a stratum. Moreover, there exists a dense stratum $R$, 
the {\em regular stratum}.

The work initiated by  Goresky and MacPherson have proved that the 
right homology to study this kind of singular spaces is the 
intersection homology (see \cite{GM1} and also \cite{M}, 
\cite{K}).
 The main result of 
this works stablishes an isomorphism between the intersection homology 
of $B$ and the basic 
intersection cohomology of $\mathcal{F}$.

\subsection{Without holonomy}
In this section, we are interested in the case where the conical 
foliation $\mathcal{F}$
is proper leaves without holonomy.
We  write $B$ the leaf space and  $\pi \colon M \to B$  the 
canonical projection.
 We shall say that a conical foliated manifold $(M,\mathcal{F})$ is a  {\em singular fibration} 
when, for each singular singular stratum $S$ of $ {\sfm{S}}_{\mathcal{F}}$, 
we have
\begin{itemize}
    
\item[a)] the restriction $\pi \colon S 
\to \pi(S)$ is a smooth fiber bundle ({\em holonomy condition for 
$\mathcal{F}_S$}),
\item[b)]  the link $(\S^n,\mathcal{G})$ of $S$ is a singular 
fibration ({\em inductive condition})
and
\item[c)] the trace of a leaf of $\mathcal{F}$ on $\S^n$ is 
connected ({\em holonomy condition for $S$ on $M$}).
\end{itemize}

\nt For this kind of foliations we have the following result. 
\begin{proposition}
    \label{noholo}
    The leaf space  of a singular fibration 
     is a stratified pseudomanifold.
    \end{proposition}   
    \pro
     Let $(M,\mathcal{F})$ be  the singular fibration and let $B$ be the leaf 
     space. We prove first that $B$ is a paracompact space. This comes 
from the fact that,
being $M$ second countable, we can find a basis $\{ (\phii,U)\}$ of 
the topology made 
up with good conical charts. The induced family $\{ (\psi,\pi(U) \}$ 
is a countable basis of $B$; the orbit space $B$ is therefore second 
countable. It is also clearly regular. This implies that $B$ is 
paracompact  (cf. \cite{Wi}). 

From a) and from the fact that the family  $ {\sfm 
S}_{\mathcal   F}$ is a locally finite then we get that the family
$\sfm{S}_{B} = \{ \pi(S)) \ / \ S \in   {\sfm 
S}_{\mathcal   F}\}$ 
 is a  
stratification.

Let us study the local structure of this stratification.
We proceed by induction
on the depth of the stratification. In the 
first step of the induction we have that the projection 
$\pi$ is a fibration and  $B$ is a stratified 
pseudomanifold whose depth is 0.  The Proposition becomes the usual 
deRham Theorem.

We consider, in the generic case,  a conical chart $(\phii,U)$ of a 
point of a stratum $S$. Since 
$(S,{\mathcal  F}_S)$ is a fibration, then  any point of $S$ admits a chart 
modelled on  $(\R^{m-n-1} ,{\mathcal H})) \equiv (\R^r, {\mathcal K}) 
\times (\R^s, {\mathcal I})$, where ${\mathcal K}$ is the foliation 
with one leaf and ${\mathcal I}$ is the foliation by points, with the 
property that each point of $\R^s$ belongs to a different leaf of 
${\mathcal  F}$. We shall say that 
$(\phii,U)$ is a {\em good conical chart} if it is modelled on 
$$ 
(\R^r 
\times \R^s \times c\S^n, {\mathcal K}\times {\mathcal I}\times 
c{\mathcal G}).
$$

\nt A good conical chart $\phii$ induces the 
stratified embedding (homeomorphism sending diffeomorphically the 
strata into the strata) $\psi \colon \R^s \times c \pi_{_{ \S^n}}( \S^n) 
\rightarrow 
B$ (cf. c)), where $\pi_{_{ \S^n}}$ is the canonical projection on $ \S^n$. 
Here, since $( \S^n , {\mathcal G})$ is a singular fibration (cf. 
b)), we can apply the induction hypothesis.

With the same technics, one shows directly that there exists a regular 
stratum $R$.   
\qed

We present now some examples.

\prgg
{\bf Examples}. 

a) We begin with a simple example. Fix $(p,q) \in \R^{2}$ and put $R\colon \C^{2} \to \C^{2}$ the 
map defined by $R(z_{1},z_{2}) = (e^{2\pi i p} \cdot z_{1},e^{2\pi i  q} \cdot 
z_{2})$. The suspension of $R$ defines on the 
manifold $M = \S^1 \times c\S^3$ a conical flow $\mathcal{F}$ without 
singularities.  
We have a singular fibration just in the case $p/q = 1$. When $p/q$ is irrational then 
the leaf space $B$ is not a stratified pseudomanifold (even 
paracompact!). In the other cases the leaf space $B$ is the 
cone over a  lens space and therefore an orbifold.

\medskip

b) Consider a smooth action of a compact Lie $G$ on a
smooth compact manifold $M$. Such an action defines a singular 
(Riemannian)
foliation $\mathcal F$ on $M$. Subsets of points of orbits of the 
same 
dimension are submanifolds of $M$ and define a stratification (cf. \cite{D}). 
On each stratum $S$, the induced foliation  $\mathcal{F}_{S}$ is 
a regular Riemannian foliation whose leaves are compact 
with finite holonomy. The space of leaves $S/\mathcal{F}_{S}$ is 
an orbifold.
If we assume that our action has just one type of orbit in each 
stratum and that the isotropy subgroups are 
connected then $(M,\mathcal{F})$ is a singular fibration.

\medskip

c)
The considerations in \cite{Mol} permits us to formulate the 
following 
generalization for compact singular riemannian foliated manifold 
$(M,\mathcal{F})$. The subsets consisting of points of 
leaves of the same dimension are smooth submanifolds and define a 
stratification of the manifold. On each stratum $S$ the foliation 
$\mathcal{F}$ defines 
a regular Riemannian foliation $\mathcal{F}_{S}$. The leaf space $S/\mathcal{F}_{S}$ is an orbifold.
Each leaf has finite holonomy. When 
all these holonomies are trivial, then $(M,\mathcal{F})$ is  a singular fibration.

\bigskip

The intersection homology $\lau{\IH}{\ast}{\per{p}}{B}$ has been 
introduced by Goresky-MacPherson  in order to study 
the stratified pseudomanifolds (see \cite{GM1} and also \cite{M}, 
\cite{K}). The main result of 
this works relates this intersection homology with the basic 
intersection cohomology of $\mathcal{F}$.
Recall that two perversities $\per{p}$ and $\per{q}$ on $B$ are {\em 
dual} if $\per{q}(\pi(S)) + \per{p}(\pi(S)) = codim_B (\pi(S)) - 2$ 
for each singular stratum $\pi(S)$. Notice that a perversity 
$\per{p}$ defines a perversity on $B$ by putting $\per{p}(\pi (S)) = 
\per{p}(S)$.

\begin{theorem}
\label{Compacto}
Let $(M,{\mathcal  F})$ be a  singular fibration. Put $B$ 
the orbit space. Then we have 
$$\lau{\IH}{\ast}{\per{q}}{\mf}  \cong 
\lau{\IH}{\ast}{\per{p}}{B},$$
where  $\per{p}$ and $\per{q}$ are dual perversities. 
\end{theorem}
\pro
See Section 4.
\qed

\subsection{\bf With holonomy.} We consider now the case where the 
holonomy is not trivial, in fact finite. We get the same relationship 
between the BIC of $\mathcal{F}$ and the intersection homology of 
$B$, but some technical difficulties appear. First of all we present an 
representative example.

\prgg{\bf Example}.
Consider $\S^{3}$ endowed with the riemannian 
foliation $\mathcal{H}$
defined by the Hopf action, that is, $z \cdot (u_{1},u_{2}) =  
(z \cdot u_{1}, z \cdot u_{2}) $. On the other hand, we consider $\R^{6} = 
c\S^{5}$ endowed with the singular riemannian foliation $c\mathcal{Q}$ defined 
by the  torus action $\T^{2} \times c\S^{5} \to c\S^5$ given by
$$
(z_{1},z_{2}) \cdot [(v_{1},v_{2},v_{3}),t] =[(z_{1} 
\cdot v_{1},z_{1}z_{2} \cdot v_{2},z_{1}z_{2}^{2} \cdot v_{3}),t].
$$
There are three different types of leaves: points, circles and tori. 
There are three different types of strata: a point (whose link is 
$(\S^{5},\mathcal{Q})$ ), three cylinders 
(whose link is $(\S^{3},\mathcal{H})$) and the regular stratum. In fact, the singular part is the wedge $\D^{2} \vee \D^{2} \vee \D^{2}$ 
endowed with the foliation induced by the natural action of $\sbat$ on 
each disk.
The finite group $\Z_{3}$ acts freely  on the product $\S^{3} \times c\S^{5}$ by
$$e^{2\pi i/3} \cdot ((u_{1},u_{2}),[(u_{1},u_{2},u_{3}),t]) =
((u_{1},e^{2\pi i/3} \cdot u_{2}),[(u_{3},u_{1},u_{2}),t]).$$
There are three different types of leaves:  circles, tori and 
three-dimensional tori. 
There are three different types of strata: $\S^{3}$, three copies of 
$\S^{3} \times \S^{1} \times ]0,1[$ 
and the regular stratum. 
Notice that this action permutes the  three four-dimensional strata .
Since this actions preserves the foliation $\mathcal{H} \times 
c\mathcal{Q}$, then we get on 
the quotient manifold 
$$M = \S^{3} \times_{\Z_{3}}  c\S^{5}$$ a singular riemannian foliation $\mathcal{F}$.
There are three different types of leaves:  circles, tori and 
three-dimensional tori. 
There are three different types of strata. The minimal stratum $S$ is the lens space $L(1,3) = \S^{3}/ \Z^3$, 
its link is $(L_{S},\mathcal{G}) = (\S^{5},\mathcal{Q})$. There is one 
four-dimensional stratum, the product $\S^{3} \times \S^{1} \times ]0,1[$ 
whose link is $(L_{S},\mathcal{G}) = (\S^{3},\mathcal{H})$. The last stratum is the regular stratum.

Put $B$ the leaf space $\mf$ and $\pi \colon M \to B$ the canonical 
projection. The image $\pi(S) = L(3,1)/\mathcal{F} = \S^{2}/\Z^{3}$ is not a manifold but an orbifold.

We already see a difference with the singular fibrations. Here, the 
strata are not necessarily manifolds but orbifolds
(see for example \cite{Mol}).
But there is another important difference related to the "links". 
The natural projection $\tau \colon M \to S$ is a foliated 
tubular neighborhood (the charts preserve the foliation and the 
conical structure) of $S$ with fiber $c\S^{5}$. The induced map
$\sigma \colon \pi(T) \to \pi(S)$ is a tubular neighborhood (the 
charts preserve the conical structure) 
os $\pi(S)$ with fiber the quotient $(\S^{5}/\mathcal{G})/\Z_{3}$. 
The "link" $L_{\pi(S)} $  of $\pi(S)$ is therefore $(\S^{5}/\mathcal{Q})/\Z_{3}$. 
In other words, 
$$
\hbox{we don't have } \ \ L_{\pi(S)} = L_{S} / \mathcal{G},
\ \ \hbox{ but } \ \ 
L_{\pi(S)} =( L_{S} / \mathcal{G})/ \Z^{3}.
$$
We can see that the trace of $\mathcal{F}$ on the link of $S$ are 
exactly the obits of $\Z^{3}\cdot \mathcal{Q}$.

\prgg {\bf Orbifolds category}. These  considerations lead us to introduce 
the following notions. A  stratified pseudorbifold
is  defined as the stratified pseudomanifold changing "manifold" by 
orbifold and "diffeomorphism" by "isomorphism between orbifolds". In 
other words, a {\em stratified pseudorbifold}
is  given by
\begin{itemize}
    \item a paracompact space $X$ 
    \item  a locally finite partition ${\sfm{S}}_{X}$, called  
    {\em stratification}, made up of connected smooth 
    orbifolds, called {\em strata}. 
    \end{itemize}
The local structure must be conical, that is, the local model is $U \times 
cL$, where $U$ is an orbifold, $L$ is a smaller compact stratified 
pseudorbifold. More 
exactly, each point of $X$ has a neighborhood which is homeomorphic to 
the local model, and this homeomorphism sends  a 
stratum into a stratum by an isomorphism. Moreover, there exists a dense stratum $R$, 
the {\em regular stratum}.

We also extend the notion of singular fibration.
 We shall say that a conical foliated manifold $(M,\mathcal{F})$ is a  
 {\em singular Seifert fibration} 
when, for each singular singular stratum $S$ of $ {\sfm{S}}_{\mathcal{F}}$, 
we have
\begin{itemize}
    
\item[a)] the restriction $\pi \colon S 
\to \pi(S)$ is a Seifert bundle,
\item[b)]  the link $(\S^n,\mathcal{G})$ of $S$ is a singular 
Seifert fibration 
and
\item[c)] the trace of a $\mathcal{F}$ on $\S^n$ is 
given by the orbits of $H \cdot \mathcal{G}$, where $H$ is a finite subgroup 
of
$$O(n+1,\mathcal{G}) = O(n+1)  \cap \Diff (\S^n,\mathcal{G}).$$
\end{itemize}

\begin{proposition}
    The leaf space  of a singular Seifert  fibration 
  is a stratified pseudorbifold.
    \end{proposition}  
\pro   Let $(M,\mathcal{F})$ be the singular fibration and  let $B$ be the leaf 
     space. We know form \cite{Mol} that the elements of the family $\sfm{S}_{B} = \{ \pi(S)) \ / \ S \in   {\sfm 
S}_{\mathcal   F}\}$ are orbifolds. Proceeding as in Proposition 
\ref{noholo} we get  that  $B$ is a paracompact space and that the 
family 
$\sfm{S}_{B}$   
is locally finite.

Let us study the local structure of this stratification.
We proceed by induction
on the depth of the stratification. In the 
first step of the induction we have that the projection 
$\pi$ is a Seifert fibration and  $B$ is a stratified 
pseudorbifold whose depth is 0.  

We consider, in the generic case,  a conical chart $(\phii,U)$ of a 
point of a stratum $S$. 
The local model of a point of $\pi(S)$ 
 is a product $(\R^p/\Gamma, c (\S^n/\mathcal{G})/H)$ where 
$\Gamma  \subset O(p)$ is a finite subgroup, $(S^n,\mathcal{G})$ is the 
link of $S$ and $H \subset O(n+1,\mathcal{G})$ is a finite subgroup.
An inductive argument gives that $\S^n/\mathcal{G}$ is a stratified 
pseudorbifold. Thus $(\S^n/\mathcal{G})/H$ is also a stratified 
pseudorbifold .\qed 
    
\prgg {\bf Singular riemannian foliations}. From \cite{BM} (see also 
\cite{SW1}) we know that each 
stratum $S$ of a singular riemannian foliation $(M,\mathcal{F})$ 
possesses a tubular neighborhood, called {\em foliated tubular 
neighborhood},  $\tau \colon T \to S$ with a 
with a foliated atlas 
$$
\{ \phii \colon (\tau^{-1}(U),\mathcal{F}) \to (U \times c\S^n,\mathcal{F} \times 
c\mathcal{G}) \}
$$
whose structural group is $O(n+1,\mathcal{F}) = \{ A \in 
O(n+1) \ / \ A \hbox{ preserves } \mathcal{G} \}$. Here 
$(\S^n,\mathcal{G})$ denotes the link of $S$. 

A singular riemannian foliation gives the main example of a stratified 
pseudomanifold.

    \bp
    The leaf space of a singular riemannian foliation with compact 
    leaves is a stratified pseudorbifold.
    \ep
    \pro  Put $(M,\mathcal{F})$ the singular riemannian foliation, 
    $B$ the leaf space and $\pi \colon M \to B$ the natural 
    projection. Since the the restriction of  $\mathcal{F}$ 
    to $S$ is a riemannian foliation  with 
    compact leaves  then the restriction $\pi \colon 
    S \to \pi (S)$ is a Seifert bundle (see \cite{Mol}). This gives a).
    
    For b) we notice that $(\S^n,\mathcal{G})$ is a singular 
    riemannian foliated manifold with $\depth \sfm{S}_{\mathcal{G}} < 
    \depth \sfm{S}_{\mathcal{F}} $. An inductive argument 
on the depth of the stratification gives b). In the 
first step of the induction we have that the projection 
$\pi$ is just a Seifert fibration.  The Proposition becomes the  
deRham Theorem for orbifolds.

The foliation  $\mathcal{F}$ does not always induces on $\S^n$ the 
foliation $\mathcal{G}$ (see the Example 3.2.1). Consider $\tau \colon 
T \to S$ a foliated tubular neighborhood whose structural group is 
$O(n+1, \mathcal{G})$. The trace of $\mathcal{F}$ 
on the generic fiber $\S^n$ is given by the orbits $H \cdot 
\mathcal{G}$ where $H = \{ A \in O(n+1, \mathcal{G}) \hbox{ 
preserving the leaves}\}$. Since $H$ preserves the regular stratum 
then, for dimensional reasons, we have $\dim H =0$. This gives c).
\qed

  As in the without-holonomy case we get that
  the basic intersection cohomology generalizes the intersection 
  cohomology of the leaf space.
  
\begin{theorem}
    \label{finito}
Let $(M,{\mathcal  F})$ be a  singular Seifert fibration. Put $B$ 
the leaf space. Then we have 
$$\lau{\IH}{\ast}{\per{q}}{\mf}  \cong 
\lau{\IH}{\ast}{\per{p}}{B},$$
where  $\per{p}$ and $\per{q}$ are dual perversities. 
\end{theorem}
\pro
The same procedure followed in the proof of the Theorem \ref{Compacto} 
reduces the problem to a local question on $B$.
Recall that the local model of $B$ is 
$$(\R^p/\Gamma, c \pi(\S^n)/H)$$
where 
$\Gamma  \subset O(p)$ is a finite subgroup, $(S^n,\mathcal{G})$ is the 
link of $S$,  $\pi \colon \S^n \to \pi(\S^n)$ denotes the canonical 
    projection associated to $\mathcal{G}$ and 
    $H \subset O(n+1,\mathcal{G})$ is a finite subgroup.
By retracting the first factor, we transform the problem to
 the proof of  the following statement:
  $$\hiru{H}{\ast}{\lau{\Omega}{\cdot}{\per{q}}{\pi(\S^n)/H}} \cong 
    \hiru{H}{\ast}{\Hom(\lau{SC}{\cdot}{\per{p}}{\pi(\S^n)/H},\R)}$$
Since $\S^n/\mathcal{G}$ is a stratified 
pseudorbifold, an inductive argument gives  

\begin{equation}
    \label{induc}
    \hiru{H}{\ast}{\lau{\Omega}{\cdot}{\per{q}}{\pi(\S^n)}} \cong 
    \hiru{H}{\ast}{\Hom(\lau{SC}{\cdot}{\per{p}}{\pi(\S^n)},\R)}.
    \end{equation}

   \nt A classic argument using the finiteness of $H$ (see for 
    example \cite{Br1}) gives
    $\hiru{H}{\ast}{\lau{\Omega}{\cdot}{\per{q}}{\pi(\S^n)/H}} 
\cong 
 \left( \hiru{H}{\ast}{\lau{\Omega}{\cdot}{\per{q}}{\pi(\S^n)}} \right)^{H}
   $ and
    $ \hiru{H}{\ast}{\Hom(\lau{SC}{\cdot}{\per{p}}{\pi(\S^n)/H},\R)}\cong 
    \left(  \hiru{H}{\ast}{\Hom(\lau{SC}{\cdot}{\per{p}}{\pi(\S^n)},\R)}\right)^{H}$,
    the subspaces of fixed points. The proof ends by noticing that the 
    isomorphism \refp{induc} is natural.
\qed

Let $\mathcal F$ be a conical foliation on $M$ with compact leaves, 
which we will call a compact conical foliation. Then on each stratum 
of $M$ the foliation induces a regular compact foliation. The results 
of  \cite{E,EMS} permit us to formulate the following Corollary. 

\begin{corollary}
Let $(M,{\mathcal  F})$ be a compact conical foliation. Let 
$\sfm{S}_{\mathcal F}$ be the partition of $M$ by subsets 
consisting of points of leaves of $\mathcal F$ of the same dimension. 
If one of the following is satisfied :
\begin{itemize}
\item[(i)] on each stratum the volume function of leaves is locally bounded;

\item[(ii)] in each stratum leaves of the induced foliation have finite 
holonomy;

\item[(iii)] the codimension of the foliation in each stratum is 2;
\end{itemize}
then the leaf space $B$, is a stratified pseudomanifold and, for any 
perversity $\per{p}$, $$\lau{\IH}{\ast}{\per{q}}{M/{\cal F}}  \cong 
\lau{\IH}{\ast}{\per{p}}{B},$$
where $\per{q}$ is the dual perversity of $\per{p}$. 
\end{corollary}

\medskip
\noindent
\pro The conditions (i),(ii) and (iii) assure that in each 
stratum the induced foliation is Riemannian. 
Therefore the leaf space $B$ endowed with the
stratification $\{ \pi(S) \  / \ S \in \sfm{S}_{\mathcal  F} \}$  is a 
stratified pseudomanifold. 
On each stratum the natural
projection is a Seifert fibration. The rest follows from Theorem 
\ref{finito}.

\section{Proof of the  Theorem \ref{Compacto}.}

We proceed in five steps.
\begin{itemize}
    
    \item We reduce the problem to a question on the leaf 
space $B$ by giving a presentation of the BIC of $(M,\mathcal{F})$ using a 
complex of differential forms living on $\pi(R)$: the complex 
$\lau{\Omega}{\ast}{\per{q}}{B}$.

\item We present the intersection homology of $B$ by using the 
complex 
$\lau{SC}{\per{p}}{\ast}{B}$   of
$\per{p}$-intersection chains.

\item Since these two complexes are not comparable by 
integration, we introduce  the 
subcomplex 
$\lau{LC}{\per{p}}{\ast}{B}$   of smooth
$\per{p}$-intersection chains.

\item
These complexes are related by the restriction $\J \colon Hom(\lau{LC}{\per{p}}{\ast}{B},\R) 
\TO 
Hom(\lau{SC}{\per{p}}{\ast}{B},\R)$ and the integration $\I \colon \lau{\Omega}{\ast}{\per{q}}{B}
\TO 
Hom(\lau{LC}{\per{p}}{\ast}{B},\R)$.

\item We prove that the operators $\I$  and $\J$ 
are quasi-isomorphisms.

\end{itemize}

\smallskip

\subsection{The complex 
$\boldsymbol{\Omega}^{\boldsymbol{\ast}}_{\boldsymbol{\per{q}}}(\boldsymbol{B})$.}

The result we are going to prove is a comparison between a cohomology 
defined on $M$ with a cohomology defined on $B$. In order to simplify 
the proof we are going to present the basic intersection cohomology 
of $(M,{\cal F})$ using differential forms on $B$. We define 

$$\lau{\Omega}{\ast}{\per{q}}{B} = \{\eta \in 
\hiru{\Om}{\ast}{\pi(R)} \ / \ \pi^{\ast}\eta \in 
\lau{\Omega}{\ast}{\per{q}}{\mf}\}.$$

Clearly, the differential operator 
$$
    \pi^{\ast} \colon 
\lau{\Omega}{\ast}{\per{q}}{B} \rightarrow 
\lau{\Omega}{\ast}{\per{q}}{\mf}
$$
is an isomorphism.

Consider $\{ V , W \}$ an open covering of $B$. Recall that there 
exists a subordinated partition of the unity made up of controlled 
functions, elements of  $\lau{\Omega}{0}{\per{0}}{B}$ (cf. \cite{S}). 
So, the covering $\{ 
\pi^{-1}(U),\pi^{-1}(V)\}$ is a basic covering and  we get from 
2.4.1 (b)
the Mayer-Vietoris 
short sequence
$$
0 \rightarrow \lau{\Omega}{\ast}{\per{q}}{B} \rightarrow 
\lau{\Omega}{\ast}{\per{q}}{V} \oplus 
\lau{\Omega}{\ast}{\per{q}}{W} \rightarrow
\lau{\Omega}{\ast}{\per{q}}{V\cap W} \rightarrow 0.
$$

The fact that a differential form of $\pi(R)$ lives on  
$\lau{\Omega}{\ast}{\per{q}}{B}$
is a  local question. But we have more than that. Notice that, any two 
points of $M$ can be related by a local foliated diffeomorphism (use 
good conical  charts). So, in order to verify that a differential 
form on $\pi(U\cap R)$ lives on 
$\lau{\Omega}{\ast}{\per{q}}{\pi(U)}$, where $(\phii,U)$ is a good  
conical chart, it suffices to look at $U$ and not at $\pi^{-1}(U)$. 
This implies that we have the 
differential isomorphism 
\begin{equation}
    \label{red}
    \phii^{\ast} \colon \lau{\Omega}{\ast}{\per{q}}{\pi(U)} 
\rightarrow 
    \lau{\Omega}{\ast}{\per{q}}{\R^s \times  \pi_{_{ \S^n}}( \S^n)},
    \end{equation}
where $\pi(U)$ is the leaf space of $((\pi^{-1}(\pi(U)), 
{\mathcal  F}_{\pi^{-1}(\pi(U))})$ and $\R^s \times  \pi_{_{ \S^n}}( 
\S^n)$ is the 
leaf space of $(\R^{m-n-1} \times   \S^n, {\mathcal H} \times 
c{\mathcal G})$.

\subsection{The complex 
$\boldsymbol{SC}^{\boldsymbol{\per{p}}}_{\boldsymbol{\ast}}(\boldsymbol{B})$.}

 For each integer $i$ we 
shall write $B_{i} \subset B$ the union of strata with dimension less 
or equal to $i$. The union of singular strata of $B$ is written 
$\Sigma_{B}$. We write $\Delta$ be a simplex and let $P$ be a prism (i.e. a finite product of 
simplices).

A  {\em liftable prism}  $\c \colon P \times \Delta \rightarrow 
B \times [0,1[^a$  is a continuous map verifying the two following 
properties:

\medskip

i) The restriction $\c \colon P \times \inte ( \Delta) \rightarrow 
\pi (R)\times 
[0,1[^a$ is smooth.

\medskip

ii)  Each $\c^{-1}(B_{i} \times [0,1[^a )$ is of the form  $P 
\times F_i$, where  $F_i$ is a face of $\Delta$.

\medskip

\nt The liftable prism $\c$ is {\em $\per{p}$-allowable} if it verifies

\medskip

iii) $\codim_{_B}\pi (S) \leq \codim_{_{\Delta}} F_{_{\dim \pi(S)}} + 
\per{p}(\pi(S))$ for each singular stratum $\pi(S)$ of $B$.

\bigskip

 A singular chain $\xi$ is an {\em $\per{p}$-intersection chain}  
when the chains $\xi$ and $\partial \xi$ are made up with 
$\per{p}$-allowable simplices. 
We shall write $\lau{SC}{\per{p}}{\ast}{B}$  the complex of 
$\per{p}$-intersection chains. The intersection homology 
$\lau{\IH}{\ast}{\per{p}}{B}$ of $B$ can be computed using the complex 
of  intersection chains $\lau{SC}{\per{p}}{\ast}{B}$ (cf. \cite{K}, 
\cite{S}).

\subsection{The complex 
$\boldsymbol{LC}^{\boldsymbol{\per{p}}}_{\boldsymbol{\ast}}(\boldsymbol{B})$.}

First of all we need stablish some results about the blow up of a 
standard simplex.

\prgg {\bf Linear blow-up}.  

Consider $\Delta = \Delta_0 * \cdots * \Delta_{k}$ a {\em 
decomposition} of  $\Delta$. We can think  $P 
\times \Delta$ as a {\em stratified prism} with  singular strata   
$\{P \times \Delta_0, P \times ((\Delta_{0} * \Delta_1) - \Delta_{0}), \ldots , P \times 
((\Delta_0 * \cdots * \Delta_{k-1}) - (\Delta_0 * \cdots * \Delta_{k-2})) 
\}$ and with a regular stratum $R_{P 
\times \Delta} = P \times (\Delta -(\Delta_0 * \cdots * \Delta_{k-1})) $. 
The {\em depth} of the stratified prism is $\depth (P \times \Delta) = k$.

 When this depth is strictly positive we can  
desingularise $P \times \Delta$ in the following way. 
The {\em linear 
blow-up} of $P \times \Delta$ is the smooth map $${\mathcal  L}_{P 
\times \Delta} \colon (P \times \per{c} \Delta_{0}) \times (\Delta_1 
* \cdots * \Delta_k) \longrightarrow P \times \Delta $$ defined by 
${\mathcal  L}_{P \times \Delta}(x,[x_0,t_0],y) = (x,t_{0}x_{0}+ 
(1-t_{0})y)$. Here $\per{c}\Delta_{0}$ denotes the closed cone 
$\frac{\Delta_{0}\times [0,1]}{\Delta_{0}\times 
\{ 0\}}$. We shall write $({P \times \Delta})^{blu} = (P \times \per{c} 
\Delta_{0}) \times (\Delta_1 * \cdots * \Delta_k)$, which is a 
stratified prism with $ \depth   
({P \times \Delta})^{blu} <  \depth  (P \times \Delta)$. The inverse 
image ${\mathcal  L}_{P \times     \Delta}^{-1}(R_{P \times \Delta})= 
R_{({P \times \Delta})^{blu}} -((P     \times \Delta_{0} \times \{ 1\}) 
\times (\Delta_1 * \cdots * \Delta_k))$  is a dense subset of 
$({P     \times \Delta})^{blu}$ and the restriction

   $$	{\mathcal  L}_{P \times \Delta} \colon {\mathcal  L}_{P \times 
    \Delta}^{-1}(R_{P \times \Delta}) \rightarrow R_{P \times 
\Delta}   $$
\noindent
is a diffeomorphism. The same properties hold for $int(P \times 
\Delta) \subset R_{P \times \Delta}$, here
\begin{equation}
	\label{diff}
     {\mathcal  L}_{P \times 
    \Delta}^{-1}(int(P \times \Delta)) = int(({P \times \Delta})^{blu}).
      \end{equation}
    Notice that, when $k =1$ and $\dim \Delta_{1} =0$, then 
${\mathcal  
L}_{P \times \Delta}$ itself  is in fact a diffeomorphism.

\prgg {\bf Smooth intersection homology}.

A liftable simplex $\c \colon \Delta \to B$
induces a natural decomposition on $\Delta$. Consider $\{ 
i_0,\ldots,i_k\}$ the family of indices verifying $F_i \not= F_{i-1}$ 
and put $\Delta_j$ the face of $\Delta$ with $F_{i_{j}} 
=F_{i_{j}-1}$.  This defines on 
$\Delta$ the $\c$-decomposition $\Delta = \Delta_0 * \cdots * 
\Delta_k$.
We have $  \c^{-1}(\Sigma_{B}\times [0,1[^a)= P \times 
(\Delta_{0} * \cdots * \Delta_{k-1})$.

\noindent 
The prism is {\em smooth} when it  also verifies the condition:

\medskip

iv) There exists a smooth map $\c \colon P \times \Delta 
\rightarrow M \times [0,1[^a$ with $\pi \rond \c = \c$.

\medskip

\noindent 
Since $\c^{-1}(\pi(R) \times [0,1[^a) = R_{P \times \Delta}$ 
then the 
    restriction $\c \colon R_{P \times \Delta} \rightarrow R 
\times [0,1[^a$ 
    is smooth and therefore it verifies a stronger condition than i), 
    namely
   
\medskip

    i)' The restriction $\c \colon R_{P \times \Delta} 
\rightarrow \pi(R) \times 
    [0,1[^a$ is smooth.

\bigskip

 A singular chain $\xi$ is a smooth {\em $\per{p}$-intersection chain}  
when the chains $\xi$ and $\partial \xi$ are made up with 
smooth $\per{p}$-allowable simplices. 
We shall write $\lau{LC}{\per{p}}{\ast}{B}$  the complex of 
smooth $\per{p}$-intersection chains. It will be shown in the next 
section that this complex also computes the intersection homology.

\subsection{The operators $\boldsymbol{\I}$ and $\boldsymbol{\J}$.}

The natural  inclusion  $\lau{LC}{\per{p}}{\ast}{B} 
\hookrightarrow \lau{SC}{\per{p}}{\ast}{B}$ induces the differential operator
$$
\J \colon \Hom (\lau{SC}{\per{p}}{\ast}{B},\R)  \TO \Hom 
(\lau{SC}{\per{p}}{\ast}{B},\R).
$$

 The difficulty to integrate a differential 
form $\omega \in \lau{\Omega}{\ast}{\per{q}}{B} $ on an  intersection 
chain $\xi$ lies on the fact that $\omega$ is defined only on the 
regular stratum of $B$ while $\xi$ is defined on  $B$. For this 
reason we need some preparatory results 

\medskip

 The linear blow up is compatible with the barycentric subdivision in 
the following way. Let $\nabla$ an element of the barycentric 
subdivision of $\Delta$, endowed with the  induced decomposition. 
That is, $\nabla = \nabla_{1} * \cdots * \nabla_{l}$ where $\{ i_{1}, 
\ldots , i_{l}  \} = \{ i \in \{ 1,\ldots,k\} \ / \ \Delta_0 * \cdots 
* \Delta_i \not= \Delta_0 * \cdots * \Delta_{i-1} \}$ and $ \Delta_0 
* \cdots * \Delta_{i_{j}} = \Delta_0 * \cdots * \Delta_{i_{j}-1} * 
\nabla_{j}$. Notice that $R_{\nabla} = R_{\Delta} \cap \nabla$. It 
has been proved in  \cite[page 220]{BHS}) that 

\begin{proposition}
\label{BHS}
Given an element $\nabla$ of the barycentric subdivision of $\Delta$ 
and $I \colon (P \times \nabla)^{blu} \hookrightarrow (P \times 
\Delta)^{blu}$ the natural inclusion, then there exists a smooth map 
$I \colon {P \times \nabla} \rightarrow {P \times \Delta}$ 
verifying $ I \rond {\mathcal  L}_{P \times \nabla} = {\mathcal  
L}_{P \times     \Delta}\rond I$. 
    \end{proposition}

\prgg{\bf Differential forms on $\boldsymbol{P \times \Delta}$}.

We shall write ${\bf\Pi}^{\ast}({P\times \Delta})$ the complex of 
{\em liftable forms}. When $ \depth  (P \times \Delta) =0$ then we put 
${\bf\Pi}^{\ast}({P\times \Delta}) =\hiru{\Om}{\ast}{P\times 
\Delta}$. In the generic case, we shall say that a differential form 
$\omega \in \hiru{\Om}{\ast}{R_{P \times \Delta}}$ is liftable if 
there exists a liftable form $\widehat{\omega} \in {\bf\Pi}^{\ast}({P 
\times \Delta}^{blu})$ with ${\mathcal  L}_{P \times \Delta}^*\omega = 
\widehat{\omega}$ on ${\mathcal  L}_{P \times \Delta}^{-1}(R_{P \times 
\Delta})$. Notice that,  when $k=1$ and $\dim \Delta_{1} =0$, the 
form  $\omega$ is defined in fact in $P \times \Delta$. By density, 
the lifting $\widehat{\omega}$ is unique. Then we have $\widehat{d\omega} = 
d\widehat{\omega}$. The complex ${\bf\Pi}^{\ast}({P\times \Delta})$ is 
thus differential. As usual, we shall write $${\displaystyle \int_{P 
\times \Delta}} \omega = {\displaystyle \int_{int(P \times \Delta) 
}}\omega, $$
which is not always well defined. But in our context

\begin{sublemma}
    \label{A}
Let $\omega$ be a liftable form, then the integral ${\displaystyle 
\int_{P \times \Delta} } \omega$ is finite. 
\end{sublemma}
\pro
We proceed by induction on the depth. When this depth is 0 then the 
result is clear. For the generic step we have  

\begin{equation}
    \label{lift}
    \int_{P \times \Delta} \omega = \int_{{\mathcal  L}_{P \times 
\Delta}^{-1}(\inte (P \times \Delta))}{\mathcal  L}_{P \times 
\Delta}^{\ast}\omega =\int_{\inte ({P \times \Delta})} 
\widehat{\omega}=\int_{{P \times \Delta}}\widehat{\omega}
\end{equation}
    since  \refp{diff} and $\inte (P \times \Delta) \subset R_{P \times 
\Delta}$. By   induction hypothesis this number is finite.
\qed

\begin{sublemma} 
\label{finite}
Let $\c \colon P \times \Delta \rightarrow B \times [0,1[^a$ be 
a smooth  prism. 
 If $\eta \in \hiru{\Om}{\ast}{\pi(R) \times [0,1[^a}$  with 
$\pi^{\ast}\eta \in {\bf\Pi}^{\ast}({M \times [0,1[^a})$ then    
$\c^*\eta$ is liftable.
\end{sublemma}

\pro Since $\c^*\eta \in \hiru{\Om}{\ast}{R_{P \times 
\Delta}}$ (cf.  i)') then it suffices to construct the lifting   
$\widehat{\c^*\eta} \in {\bf\Pi}^{\ast}({{P \times \Delta}})$. 
We proceed in several steps.

\bigskip

\noindent 
I - {\em Localizing $M$ and $B$}. 
Remark that for any element $\nabla$ of the barycentric subdivision 
of $\Delta$, 
the restriction $\c' \colon P \times \nabla \rightarrow B \times 
[0,1[^a$ is a smooth  prism. The statement becomes a local one. So, 
we can identify $$(M,{\mathcal  F}) \equiv ( \R^r \times  \R^s 
\times  c \S^n , {\mathcal K} \times {\mathcal I} \times c{\mathcal 
G} ), \hspace{1cm} B \equiv \R^s \times  c\pi_{_{ \S^n}}( \S^n)$$ and 
suppose that  $\Im \c$ meets $\R^{blu} \times \{ \vartheta\} \times 
[0,1[^a$. Then $\c^{-1}(\R^s \times \{ \vartheta \}\times 
[0,1[^a)  = P \times \Delta_0$. 

 Notice that a neighborhood of $\Delta_0$ on $\Delta$ is a product 
of $\Delta_0 \times c\nabla$, where $\nabla$ is a simplex. From Lemma 
\ref{lema} we get a commutative diagram

$$ 
\begin{picture}(200,70)(00,-8)

    \put(-148,27){\makebox(0,0){$ (\theequation)$}}
        \label{diag}
\put(0,0){\makebox(0,0){$ P \times \Delta$}}
\put(180,0){\makebox(0,0){$\R^{m-n-1} \times c  \S^n \times [0,1[^a$}}
\put(180,50){\makebox(0,0){$\R^{m-n-1} \times   \S^n \times 
[0,1[ \times [0,1[^a$}}
\put(0,50){\makebox(0,0){${P \times \Delta}$}}

\put(0,40){\vector(0,-1){26}} 
\put(180,40){\vector(0,-1){26}} 
\put(35,0){\vector(1,0){66}}
\put(35,50){\vector(1,0){60}}

\put(68,8){\makebox(0,0){$\c$}}
\put(65,61){\makebox(0,0){$\widehat{ \c}$}}
\put(-20,27){\makebox(0,0){${\mathcal  L}_{P \times \Delta}$}} 
\put(220,27){\makebox(0,0){$P\times \hbox{ Id }_{[0,1[^a}$}} 
\end{picture} 
$$ 
where  $\widehat\c$ is smooth. 

\bigskip

\noindent II - {\em Blowing up $\c$}. Consider now the 
continuous map 
$$
\widehat\c =\widehat\pi \rond  \widehat\c  \colon {P \times 
\Delta}
\rightarrow \R^s \times  \pi_{_{ \S^n}}( \S^n) \times [0,1[^{a+1},
$$ 
where $\widehat\pi \colon \R^{m-n-1} \times   \S^n \times 
[0,1[^{a+1}\rightarrow \R^s \times  \pi_{_{ \S^n}}( \S^n) \times 
[0,1[^{a+1}$ is the projection defined by 

$$
\widehat\pi (x_{1}, \ldots , x_{r}, y_{1}, \ldots, y_{s}, \theta, 
t_0,t_{1}, \ldots, t_{a}) = ( y_{1}, \ldots, y_{s}, \pi_{_{ 
\S^n}}(\theta), 
t_0,t_{1}, \ldots, t_{a}).
$$
Let us see that $\widehat\c$ is a smooth prism.

\medskip

i) Since $\widehat{\pi}$ and $\widehat\c$ are smooth it suffices 
to prove 
that
$$
\widehat\c^{-1}(\R^s \times  \pi_{_{ \S^n}}( \Sigma_{{\mathcal 
G}}) \times 
[0,1[^{a+1}) \subset
((P \times \per{c}\Delta_{0}) \times ((\Delta_{1} * \cdots * 
\Delta_{k}) -\inte (\Delta_{1} * \cdots * \Delta_{k})).
$$

This comes from
$$
\widehat\c^{-1}\widehat{\pi}^{-1}(\R^s \times  \pi_{_{ 
\S^n}}(\Sigma_{{\mathcal G}}) \times [0,1[^{a+1})
=
\widehat\c^{-1}(\R^{m-n-1} \times  \Sigma_{{\mathcal G}} \times 
[0,1[^{a+1})
= 
$$
$$
\widehat\c^{-1}(P\times \hbox{ Id }_{[0,1[^a})^{-1}(\R^{m-n-1} 
\times  
c\Sigma_{{\mathcal G}} \times [0,1[^{a})
=
{\mathcal  L}_{P \times \Delta}^{-1}\c^{-1}(\R^{m-n-1} \times  
c\Sigma_{{\mathcal G}} \times [0,1[^{a})
=
$$
$$
{\mathcal  L}_{P \times \Delta}^{-1}\c^{-1}\pi^{-1}(\R^s \times  
c\pi_{\S^n}(\Sigma_{{\mathcal G}}) \times [0,1[^{a})
=
{\mathcal  L}_{P \times \Delta}^{-1}\c^{-1}(\R^s \times  
c\pi_{\S^n}(\Sigma_{{\mathcal G}}) \times [0,1[^{a})=
$$
$$
{\mathcal  L}_{P \times \Delta}^{-1}(P  \times ( \Delta_{0}* 
\cdots * \Delta_{k-1})  \subset 
((P \times \per{c}\Delta_{0}) \times ((\Delta_{1} * \cdots * 
\Delta_{k}) -\inte (\Delta_{1} * \cdots * 
\Delta_{k})).
$$

\medskip

ii) Proceeding as before we get 

$$\widehat\c^{-1}(\R^s \times  ((\pi_{_{ \S^n}}( 
\S^n))_{j}\times [0,1[^{a+1}  ) = {\mathcal  L}_{P \times 
\Delta}^{-1}\c^{-1}(\R^s \times 
(c\pi_{_{ \S^n}}( \S^n))_{j+1} \times [0,1[^{a}  )=
$$
$$
{\mathcal  L}_{P \times \Delta}^{-1}(P \times (\Delta_{0} *  \cdots * 
\Delta_{h}))=
{P \times (\Delta_{0} *  \cdots * \Delta_{h})},
$$
for some $h \in \{ 1, \ldots, k\}$.

\medskip

iii) By construction.

\bigskip

\noindent 
III - {\em Lifting $\eta$}. The differential forms $\pi^{\ast}\eta 
\in {\bf\Pi}^{\ast}({\R^{m-n-1} \times c  \S^n \times [0,1[^a})$ 
lifts into the differential form $\widehat{\pi^{\ast}\eta} \in 
{\bf\Pi}^{\ast}({\R^{m-n-1} \times   \S^n \times [0,1[^{a+1}})$. 
Since $\pi^{\ast}\eta$ is basic and the restriction of $P\times 
\hbox{ Id }_{[0,1[^a}$ to $\R^{m-n-1} \times R_{0} \times ]0,1[ 
\times [0,1[^a$, where $R_0$ is the regular stratum of $ \S^{m-n-1}$, 
is the identity then $\widehat{\pi^{\ast}\eta}$ is also basic. There 
exists $\widehat\eta \in 
\hiru{\Om}{\ast}{\R^s \times \pi_{_{ \S^n}}(R_{0}) \times 
[0,1[^{a+1}}$ with 
$\widehat\pi^{\ast}\widehat\eta = \widehat{\pi^{\ast}\eta}$. The 
differential form $\widehat\eta$ is in the conditions of the Lemma.

\medskip

\noindent 
IV - {\em Final step}. 
We proceed by induction on $ \depth  {\sfm{S}}_{{\mathcal  F}}$. The 
result is clear when this depth is 0, that is, when $B=\pi(R)$. For 
the generic case, notice 
that $ \depth  {\sfm{S}}_{{\mathcal H} \times {\mathcal G}} <   \depth  
{\sfm{S}}_{{\mathcal H} \times c{\mathcal G}}$ then the induction 
argument gives that $\widehat\c^{\ast}\widehat\eta$ is liftable. 
It remains to prove that $\widehat{\c^{\ast}\eta} = \widehat{\bf 
c}^{\ast}\widehat\eta$, that is, ${\mathcal  L}_{P \times 
\Delta}^{\ast} \c^{\ast}\eta = \widehat{\bf 
c}^{\ast}\widehat\eta$ on ${\mathcal  L}_{P \times \Delta}^{-1} 
(R_{P\times \Delta})$.
Since $\c\rond {\mathcal  L}_{P \times \Delta} = (P\times \hbox{ 
Id }_{[0,1[^a}) \rond \widehat\c$ then we have

$$
{\mathcal  L}_{P \times \Delta}^{\ast} \c^{\ast}\eta = 
{\mathcal  L}_{P \times \Delta}^{\ast}\c^{\ast}\pi^{\ast}\eta= 
\widehat\c^* (P\times \hbox{ Id }_{[0,1[^a})^*\pi^{\ast}\eta
$$

on ${\mathcal  L}_{P \times \Delta}^{-1}(R_{P \times \Delta})$. Since 
$\pi^*\eta$ is liftable then $$(P\times \hbox{ Id 
}_{[0,1[^a})^{\ast}\pi^{\ast}\eta =\widehat{\pi^*\eta}= 
\widehat\pi^{\ast}\widehat\eta$$ on 
$\R^{m-n-1} \times ( \S^n - \Sigma_{ \S^n}) \times ]0,1[ \times 
[0,1[^a$,
So

$$
{\mathcal  L}_{P \times \Delta}^{\ast} \c^{\ast}\eta = 
\widehat\c^{\ast}\widehat \pi^{\ast}\widehat\eta = \widehat{\bf 
c}^*\widehat\eta,
$$ 

on $\widehat\c^{-1}(\R^{m-n-1} \times ( \S^n - \Sigma_{ \S^n}) 
\times ]0,1[ \times 
[0,1[^a) = {\mathcal  L}_{P \times \Delta}^{-1}(R_{P \times \Delta})$.
\qed

\bigskip

Given  a form $ \eta 
\in 
 \lau{\Omega}{\ast}{\per{p}}{B}$ and a  smooth $\per{p}$-allowable 
simplex  $c \colon \Delta \rightarrow B$ we can define  the integral

 $$\I(\omega)(c)= {\displaystyle \int_{\Delta} c^{\ast} \eta}$$ 
(cf. Lemma  \ref{A} and Lemma \ref{finite}). 
It remains to prove that the operator
$$\I  \colon 
\lau{\Omega}{\ast}{\per{q}}{B} \rightarrow   
Hom(\lau{LC}{\per{p}}{\ast}{B},\R).$$
is a differential one. We also need some preparatory results.

\prgg {\bf Boundary}. 

There are two types of (one codimensional) 
faces on $({P \times \Delta})^{blu}$.

\begin{itemize}
    \item[T1)] The faces $({Q \times  \nabla})^{blu} $, where $Q$ is a 
face of $P$ and $\nabla =  \Delta$ or $Q=P$ and $\nabla$ is a face of 
$\Delta$. The restriction of ${\mathcal  L}_{P \times \Delta}$ is the 
linear blow up ${\mathcal  L}_{Q \times \nabla}$.
    
    \item[T2)] The face $F_{P\times \Delta} = (P\times \Delta_{0} 
\times \{ 1\} ) \times (\Delta_1 * \cdots * \Delta_k)$. The 
restriction of ${\mathcal  L}_{P \times \Delta}$  is just the 
canonical projection over $ P \times \Delta_{0}$.
\end{itemize}

\noindent    
The faces of type T1)  run over the boundary of $P \times \Delta$. 
The face $F_{P\times \Delta}$ is the extra face produced by the 
linear blow up.  We have the formula
\begin{equation}
\label{bound}
\partial (({P\times \Delta})^{blu})=({\partial(P \times \Delta)})^{blu} 
+  F_{P\times \Delta}.
\end{equation}

\prgg{\bf More differential forms on $\boldsymbol{P \times \Delta}$}.
We shall write ${\Gamma}^{\ast}{(P\times \Delta)}$ the complex of 
{\em regular forms}. When $ \depth  (P \times \Delta) =0$ then we put 
${\Gamma}^{\ast}({P\times \Delta}) =\hiru{\Om}{\ast}{P\times 
\Delta}$. In the generic case, we shall say that a liftable form 
$\omega \in {\bf\Pi}^{\ast}({P \times \Delta})$ is regular  if  
$\widehat{\omega} \in {\bf\Gamma}^{\ast}({{P \times \Delta}})$ and if its 
restriction to $\inte F_{P \times \Delta}$ vanishes. Notice that 
${\bf\Gamma}^{\ast}{(P\times \Delta)}$ 
is a differential complex. From Proposition \ref{BHS} we know that 
the restriction of a regular form to a element of the barycentric 
subdivision is again regular. For these kind of forms we have the 
following Stoke's Theorem.

\begin{sublemma}
\label{regular}
Let $\omega$ a regular form, then  ${\displaystyle \int_{P \times 
\Delta} } d \omega= {\displaystyle \int_{\partial(P \times \Delta)} 
}  \omega$.
\end{sublemma}
\pro
We proceed by induction on the depth. When this depth is 0 then the 
result is clear. For the generic step we notice that $\omega$ is 
defined in the interior of  
any face $F$ of  $P \times \Delta$ except in the case where $\dim 
\Delta_{k} = 0$ and $F = \Delta_{0} * \cdots * \Delta_{k-1}$. But 
here $\omega$ extends to  a form on $\inte F$. To see that, we apply 
the Proposition \ref{BHS} to reduce the problem to $\Delta = 
\Delta_{0} * \Delta_{1}$ and $\dim \Delta_{1}=0$. Here we 
know that ${\mathcal  L}_{P \times \Delta}$ is a diffeomorphism and 
that $\widehat{\omega}$ is defined everywhere, then $\omega$ is defined 
everywhere. The two terms of the equality to show make sense.

The induction hypothesis gives 

$$\int_{{P \times \Delta}}\widehat{d\omega} = \int_{\partial({P \times 
\Delta})}\widehat{\omega}$$

 and therefore (cf. \refp{bound} and \refp{lift})

 $$
\int_{P \times \Delta }d\omega = \int_{{P \times 
\Delta}}\widehat{d\omega} = 
\int_{\partial({P \times \Delta})}\widehat{\omega} = \int_{{\partial (P 
\times 
\Delta)} }\widehat{\omega} + \int_{ F_{P\times \Delta}}\widehat{\omega} = 
\int_{\partial (P \times \Delta) }\omega + \int_{ F_{P\times 
\Delta}}\widehat{\omega}.
 $$
Now it suffices to notice that $\widehat{\omega}$ vanishes on $F_{P 
\times \Delta}$.
 \qed

\begin{sublemma} 
\label{reu}
Let $\c \colon P \times \Delta \rightarrow B \times [0,1[^a$ be 
a smooth $\per{p}$-allowable prism. If $\eta \in 
\lau{\Omega}{\ast}{\per{q}}{B}$ then   $\c^*\eta$ is regular. 
\end{sublemma}

\pro
Since $\c^*\eta$ is liftable, it remains to prove that 
$\widehat{\c^*\eta}$ 
vanishes on $\inte F_{P\times \Delta}$. In fact, it suffices to prove 
that 
 that,  we have
 \begin{equation}
   \label{perve}
\widehat\c^*\widehat\eta(u_1, \ldots, u_{a}, 
v_{1},\ldots,v_{b},0,w_{1}, \ldots 
w_{c})=0,
 \end{equation}

where   $\{ u_{1},\ldots , u_{a = \dim P} \}$ are tangent vectors to 
$\inte P$,  
$\{ v_{1},\ldots , v_{b = \dim \Delta_{0}} \}$ are tangent vectors to 
$\inte \Delta_{0}$ and $\{ w_{1},\ldots , w_{c = \dim \Delta_{1} * 
\cdots *\Delta_{k}} \}$ are tangent vectors to $\inte (\Delta_{1} * 
\cdots *\Delta_{k})$.

Since the question is finally local, we can proceed as before and 
identify $(M,{\mathcal  F})$ with $ ( \R^r \times  \R^s \times  c 
\S^n, {\mathcal K} \times {\mathcal I} \times c{\mathcal G} )$, 
identify $B$ with  $ \R^s \times  c\pi_{_{ \S^n}}( \S^n)$, identify 
$S$ with $\R^{m-n-1} \times \{ \vartheta \}$ and suppose that $\Im 
\c$ meets $\R^{blu} \times \{ \vartheta\} \times [0,1[^a$. Then 
$\c^{-1}(\R^s \times \{ \vartheta \}\times [0,1[^a)  = P \times 
\Delta_0$.
The diagram \refp{diag} becomes 

$$ 
\begin{picture}(200,70)(30,-2)

\put(55,0){\makebox(0,0){$ P \times \Delta_{0}$}}
\put(240,0){\makebox(0,0){$\R^{m-n-1} \times \{ \vartheta \} \times 
[0,1[^a$}}
\put(240,50){\makebox(0,0){$\R^{m-n-1} \times   \S^n \times 
\{ 0\}  \times [0,1[^a$}}
\put(55,50){\makebox(0,0){$F_{P \times \Delta}$}}

\put(55,40){\vector(0,-1){26}} 
\put(240,40){\vector(0,-1){26}} 
\put(95,0){\vector(1,0){66}}
\put(95,50){\vector(1,0){60}}

\put(128,8){\makebox(0,0){$\c$}}
\put(125,61){\makebox(0,0){$\widehat{ \c}$}} 
\put(10,27){\makebox(0,0){${\mathcal  L}_{P \times \Delta}$}} 
\put(290,27){\makebox(0,0){$P\times \hbox{ Id }_{[0,1[^a}$}} 
\end{picture} 
$$ 
\noindent
Here ${\mathcal  L}_{P \times \Delta}(x_1, \ldots, x_{a}, 
y_{1},\ldots,y_{b},0,z_{1}, \ldots z_{c}) = (x_1, \ldots, x_{a}, 
y_{1},\ldots,y_{b})$. The prism $\widehat\c$ verifies i)' and 
therefore $\widehat\c$ sends $\inte F_{P \times \Delta}$ into 
$\R^{m-n-1} \times  ( \S^n - \Sigma_{{\mathcal G}}) \times \{ 0\}  
\times [0,1[^a$. Since $\widehat\c = \widehat\pi \rond 
\widehat\c$  and $\widehat\pi^{\ast}\widehat\eta = 
\widehat{\pi^{\ast}\eta}$ then the condition \refp{perve} is 
equivalent to 
$$
\widehat{\pi^{\ast}\eta}(\widehat\c_*u_1, \ldots, \widehat{\bf 
c}_*u_{a}, 
\widehat\c_*v_{1},\ldots,\widehat\c_*v_{b},0,\widehat{\bf 
c}_*w_{1}, \ldots 
\widehat\c_*w_{c})=0.
$$
Since $\codim_{_{\Delta}} F_{\dim \pi(S)}=\codim_{_{\Delta}} F_{s} = 
\codim_{_{\Delta}}\Delta_{0} = c+1$ then the condition iv) implies that 

$$
c \geq \codim_{_B} \pi (S) - \per{p}(\pi(S)) -1 = \per{q}(\pi(S)) + 1 
> 
\per{q}(\pi(S)).
$$

Finally, since $({\mathcal  L}_{P \times \Delta})_{\ast}w_{j}=0$ then 
$\widehat\c_{\ast}w_{j}$ is a vector of $\R^{m-n-1} \times ( 
\S^n - \Sigma_{{\mathcal G}}) \times \{ 0\} \times [0,1[^a$ tangent 
to the fibers of $(P \times
 \hbox{ Identity }_{[0,1[^a})$ and therefore we get  (10).
  \qed
 
  \bigskip
 The Lemma \ref{regular} and the Lemma \ref{reu} show that the operator
$$\I  \colon 
\lau{\Omega}{\ast}{\per{q}}{B} \rightarrow   
Hom(\lau{LC}{\per{p}}{\ast}{B},\R)$$
is a differential operator.

\subsection{The quasi-isomorphisms $\boldsymbol{\I}$ and 
$\boldsymbol{\J}$.}

Consider the statement  $Q(B)$ about the leaf space of 
singular fibrations:
 $$ Q(B) = "\lau{\Omega}{\ast}{\per{q}}{B} \stackrel{\I}{\rightarrow} 
Hom(\lau{LC}{\per{p}}{\ast}{B},\R) \stackrel{\J}{\longleftarrow} 
Hom(\lau{SC}{\per{p}}{\ast}{B},\R)\hbox{ are quasi-isomorphisms}",$$
and we shall prove it by using the Bredon's Trick \cite[pag. 289]{Br}. 

\bigskip

\nt {\bf Bredon's Trick}. {\em 
    Let $X$ be a paracompact topological espace
and  let $\{ U_\alpha\}$ be an open covering, closed for 
finite intersection. Suppose that $Q(U)$ is a statement about open subsets 
of $M$, satisying the following three properties:

\medskip

(BT1)  \ $Q(U_\alpha)$ is true for each $\alpha$;

\smallskip

(BT2) \ $Q(U)$, $Q(V)$ and $Q(U\cap V )$ $\Longrightarrow$ $Q(U\cup V 
)$, where $U$ and $V$ are open subsets of $M$; 

\smallskip

(BT3) \ $Q(U_i) \Longrightarrow Q({\displaystyle \bigcup_i }U_i)$, 
where $\{ U_i\}$ is a disjoint family of open subsets of $M$.

\medskip

\nt Then $Q(X)$ is true.
}
\bigskip

We proceed by 
induction on $ \depth  {\sfm{S}}_{\mathcal   F}$, then this depth is 0 
then $Q(B)$ is the usual de Rham Theorem. Let us suppose that $Q(B)$ 
is proved when $ \depth  {\sfm 
S}_{\mathcal   F} < \ell$. We proceed in several steps.
\bigskip

\noindent
{\bf (i)} {\em The foliated manifold $(M,{\mathcal  F})$ is 
$(\R^u \times \R^v \times c  S^w,{\mathcal K} \times {\mathcal I} 
\times c {\mathcal G})$ where ${\mathcal K} $  is the foliation with 
one leaf,  ${\mathcal I}$ is the foliation by points and  $( S^v, 
{\mathcal G})$ is a singular fibration with $\depth  {\cal 
S}_{\mathcal   G}< \ell$}.
 
Using the canonical projection $\Pr \colon \R^v \times c\pi_{_{ S^w}}( 
S^w )
 \rightarrow c\pi_{_{ S^w}}( S^w ) $, restricted to the regular 
part,  we already know that 
$$\bi{\Pr}{\ast} \colon \lau{\Omega}{\ast}{\per{q}}{c\pi_{_{ S^w}}( S^w) 
} \rightarrow \lau{\Omega}{\ast}{\per{q}}{\R^v \times c\pi_{_{ S^w}}( 
S^w) }$$ 
is a quasi-isomorphism (cf. Proposition \ref{local4} and (8)).
Following the same procedure as in \cite{BHS} we prove that the  
operators 

$$\bi{\Pr}{\ast} \colon Hom(\lau{SC}{\per{p}}{\ast}{c\pi_{_{S^w}}( 
S^w)},\R)\rightarrow 
Hom(\lau{SC}{\per{p}}{\ast}{\R^v \times c\pi_{_{ S^w}}( S^w)},\R)$$

\nt and

$$\bi{\Pr}{\ast} \colon Hom(\lau{LC}{\per{p}}{\ast}{c\pi_{_{ S^w}}( 
S^w)},\R) \rightarrow Hom(\lau{LC}{\per{p}}{\ast}{\R^v \times c\pi_{_{ 
S^w}}( S^w)},\R)$$
 \noindent
are quasi-isomorphisms. Now, since $\bi{\Pr}{\ast}$ commutes with $\I$ and 
$\J$, 
then $Q(c\pi_{_{ S^w}}( S^w) ) \Longrightarrow Q(\R^v \times c\pi_{_{ 
S^w}}( S^w) )$. It remains to prove $Q(c\pi_{_{ S^w}}( S^w) )$.

From Proposition \ref{local5} and (8) we know that the canonical 
projection $\Pr \colon  ( S^w -\Sigma_{{\mathcal G}}) \times ]0,1[ 
\rightarrow 
( \S^n -\Sigma_{{\mathcal G}})$ induces the quasi-isomorphism

$$
\lau{\IH}{i}{\per{q}}{c\pi_{_{ S^w}}( S^w )} = \left\{
\begin{array}{cl}
\lau{\IH}{i}{\per{p}}{\pi_{_{ S^w}}( S^w )} & \hbox{if } i \leq \per{q}(\{ 
\vartheta \}) \\
0 & \hbox{if } i > \per{q}(\{ \vartheta \}) 
\end{array}
\right .
$$

Following the same procedure as in \cite{BHS} we prove that 
$\bi{\Pr}{\ast}$ induces  the isomorphisms 

$$H^{i}({\lau{LC}{\per{p}}{\ast}{c\pi_{_{ S^w}}( S^w)}} =
\left\{
\begin{array}{cl}
H^{i}({\lau{LC}{\per{p}}{\ast}{\pi_{_{ S^w}}( S^w)}} & \hbox{if 
} i \leq 
k -1 - \per{p}(\{ \vartheta \}) \\
0 & \hbox{if } i > k -1 -\per{p}(\{ \vartheta \}) 
\end{array}
\right.
$$
and 
$$
H^{i}({\lau{SC}{\per{p}}{\ast}{c\pi_{_{ S^w}}( S^w)}} =
\left\{
\begin{array}{cl}
{H}^{i}({\lau{SC}{\per{p}}{\ast}{\pi_{_{ S^w}}( S^w)}}& \hbox{if 
} i \leq 
k -1 - \per{p}(\{ \vartheta \}) \\
0 & \hbox{if } i > k -1 -\per{p}(\{ \vartheta \}) ,
\end{array}
\right.
$$
 where $k=\dim \pi_{_{ S^w}}( S^w)$. Now, since $\I$ and $\J$ commute  
with $\bi{\Pr}{\ast}$ it suffices to apply  $Q(\pi_{_{ S^w}}( S^w ))$ and 
to take into account that the  perversities $\per{p}$ and $\per{q}$ 
on $B$ are dual.

\bigskip

\noindent
{\bf (ii)} {\em The foliated manifold $(M,{\mathcal  F})$ is an open 
subset  of  $(\R^u \times \R^v \times c  S^w,{\mathcal K} \times 
{\mathcal I} \times c {\mathcal G})$ where ${\mathcal K} $  is the 
foliation with one leaf,  ${\mathcal I}$ is the foliation by points 
and  $( S^v, {\mathcal G})$ is a singular fibration with $ \depth  
{\sfm{S}}_{\mathcal   G}< \ell$}.

Since $B$ is paracompact we can apply the Bredon's Trick to the 
following 
covering
$$
\{ V \times c_\epsilon \pi_{_{ S^w}}( S^w)  \ / \ V = ]a_1,b_1[\times 
\cdots \times ]a_{v},b_{v}[ \subset \R^v, \epsilon \in [0,1[ \}
 \cup \{ \pi(U)  \ / \ U  \subset \R^{u+v}   \times S^w \times ]0,1[ \hbox{ open }\},
 $$
\noindent
 where  $c_\epsilon \pi_{_{ S^w}}( S^w) = \pi_{_{ S^w}}( S^w) \times 
[0,\epsilon [   \ / \ \pi_{_{ S^w}}( S^w) \times\{ 0 \}$. This family 
is closed for finite intersections. Let us verify the three 
conditions (BT1-3).

\medskip
(BT1) Apply (i)
to $V \times c_\epsilon \pi_{_{ S^w}}( S^w)$  and apply the induction hypothesis to 
$\pi(U)$ noticing that$ \depth  {\sfm{S}}_{{\mathcal  F}_{U}} <  \ell$. 

 \medskip
 
 (BT2) This is Mayer-Vietoris.
 
 \medskip
 
 (BT3) By construction.
 
 \bigskip
 
\noindent
{\bf (iii)}{\em The depth of the foliated manifold $(M,{\mathcal  
F})$ is $\ell$.} 

Since $B$ is paracompact we can apply the Bredon's 
Trick relatively to the following covering:
$$
\{ V \subset B  \ / \ V \hbox{ is an open subset of $\pi(U)$ where 
$(\phii,U)$ is a good conical chart} \}.  $$
 This family is closed for finite intersections. Let us verify the 
three conditions (BT1-3).

\medskip

(BT1) Apply\refp{red}  and (ii)
using the fact that $ \ell = \depth \sfm{S}_{{\mathcal K} \times 
{\mathcal I} \times c {\mathcal G}} = \depth \sfm{S}_{ {\mathcal G}}+1$.

 \medskip
 
 (BT2) This is Mayer-Vietoris.
 
 \medskip
 
 (BT3) By construction.

\end{document}